%% file: HuangEtAl_Expensive_control_SOS_archive_2016.tex
\documentclass{article}

\pdfoutput=1

\input{Eedefs}
\usepackage{graphics} 
\usepackage{epsfig} 
\usepackage{mathptmx} 
\usepackage{times} 
\usepackage{amsmath} 
\usepackage{amssymb}  
\usepackage{bm}
\usepackage{graphicx}
\usepackage{float}
\usepackage{color}
\usepackage{framed} 
\definecolor{shadecolor}{rgb}{0.8,0.8,0.8}
\usepackage[justification=centering]{caption}
\usepackage{epstopdf}
\usepackage{array}

\begin{document}
\newcommand{\mymtrx}[1]{\mathrm{#1}}
\newcommand{\todo}[1]{
\begin{shaded}
\setlength{\columnwidth}{1in}
{#1}
\end{shaded}
}
\newcommand{\commentout}[1]{}
\newcommand{\bfG}{{\mathbf G}}
\newcommand{\bfH}{{\mathbf H}}

\title{Expensive control of long-time averages using sum of squares and its application to a laminar wake flow}

\author{Deqing Huang$^1,$  Bo Jin$^2,$ Davide Lasagna$^3,$ Sergei Chernyshenko$^2$\footnote{Corresponding author, s.chernyshenko@imperial.ac.uk},\\ and Owen Tutty$^3$}

\maketitle
\begin{center}
{\footnotesize
$^1$School of Electrical Engineering, Southwest Jiaotong University, Chengdu, 610031, China

$^2$Department of Aeronautics, Imperial College London, SW7 2AZ, UK

$^3$Engineering and the Environment, University of Southampton, 
 SO17 1BJ, UK%
 
}
\end{center}

\begin{abstract}
The paper presents a nonlinear state-feedback control design approach for long-time average cost control, where the control effort is assumed to be expensive. The approach is based on sum-of-squares and semi-definite programming techniques. It is applicable to dynamical systems whose right-hand side is a polynomial function in the state variables and the controls. The key idea, first described but not implemented in (Chernyshenko \emph{et al.} Phil. Trans. R. Soc. A, 372, 2014), is that the difficult problem of optimizing a cost function involving long-time averages is replaced by an optimization of the upper bound of the same average. As such, controller design requires the simultaneous optimization of both the control law and a tunable function, similar to a Lyapunov function. The present paper introduces a method resolving the well-known inherent non-convexity of this kind of optimization. The method is based on the formal assumption that the control is expensive, from which it follows that the optimal control is small. The resulting asymptotic optimization problems are convex.  The derivation of all the polynomial coefficients in the controller is given in terms of the solvability conditions of state-dependent linear and bilinear inequalities. The proposed approach is applied to the problem of designing a full-information feedback controller that mitigates vortex shedding in the wake of a circular cylinder in the laminar regime via rotary oscillations. Control results on a reduced-order model of the actuated wake and in direct numerical simulation are reported.
\end{abstract}




\section{Introduction}
Global stabilization of dynamical systems is of importance in system theory and engineering~\cite{Kh:02,An:02}, but it is sometimes difficult or impossible to synthesize a global stabilizing controller~\cite{Va:02}. The reasons could be the poor controllability of the system~\cite{Di:11,Gu:13,Gu:14}, the input/output constraints in practice~\cite{Bl:99}, time delay~\cite{Sun:11,Sun:13}, and/or the involved large disturbances~\cite{Ki:06}, etc.  Moreover, in many applications, full stabilization, while possible, carries too high penalty due to the cost of the control.

Instead, reducing the long-time average of the cost, even by a small amount, might be more realistic, especially when the control is expensive~\cite{Dan:12}. Long-time-average cost analysis and control is often considered in irrigation, flood control, navigation, water supply, hydroelectric power, computer communication networks, and many other cases~\cite{Du:10,Bo:92}. Systems including stochastic factors are often controlled in the sense of long-time averages. In~\cite{Ro:83}, a summary of long-time-average cost problems for continuous-time Markov processes is given. In~\cite{Me:00}, the long-time-average control of a class of problems that arises in the modelling of semi-active suspension systems was considered. However, in certain cases the computational complexity of direct calculation of converged averages may be prohibitive.
To overcome this difficulty, we adopt the perspective first described in~\cite{Ph:14}, where instead of considering the long-time average as the cost for system analysis and control design, we use as the cost the upper bound of the long-time average. For bounds tight enough the control reducing the upper bound will also reduce the long-time average itself. In this paper, we describe a numerically tractable approach for long-time average cost control (LTACC) for a class of nonlinear dynamical systems with the right-hand side given by polynomials in the state variables and in the control inputs. For such systems, the sum-of-squares (SOS) decomposition of polynomials and semidefinite programming (SDP) techniques are used to estimate and optimize bounds on long-time averages.

The SOS methods apply to systems defined by a polynomial vector field. Recent results on SOS decomposition have transformed the verification of non-negativity of polynomials into SDP, hence providing promising algorithmic procedures for stability analysis of polynomial systems. However, using SOS techniques for optimal control, as for example in~\cite{Pr:02,Zh:07,Ma:10}, is subject to a generic difficulty: while the problem of optimizing the candidate Lyapunov function certifying the stability for a closed-loop system for a given controller and the problem of optimizing the controller for a given candidate Lyapunov function are reducible to an SDP and thus, are tractable, the problem of simultaneously  optimizing both the control and the Lyapunov function is non-convex. Iterative procedures were proposed to overcome this difficulty~\cite{Zh:07,Zh:09,Ng:11},
where the performance of optimization is highly dependent on the initial guess for the tuning variables.



While designing a controller that reduces an upper bound does not involve a Lyapunov function, it does involve a similar tunable function, and it shares the same difficulty of non-convexity. In the present work this difficulty is overcome by making use of the idea of expensive control~\cite{Chen:07}. We propose a polynomial-type state-feedback controller design scheme for the long-time average upper-bound control, where there is a small-amplitude parameter characterizing the expensive controller. Expanding the tunable function and the bound in the small parameter leads to convex problems. The derivation of all the polynomial coefficients in controller is given in terms of the solvability conditions of state-dependent linear and bilinear inequalities. Notice the significant conceptual difference between our approach and the studies of control by small perturbations, often referred to as tiny feedback, see for example~\cite{tc:93}.

As an illustrative example, the proposed control design scheme is applied to the problem of mitigating developed vortex shedding in the two-dimensional incompressible flow past a circular cylinder in the laminar regime at Reynolds number equal to 100. The flow is controlled by  rotary motions of the cylinder. This configuration is often used as a benchmark problem to test modeling/control algorithms and strategies for fluid flows. The control perspective is to reduce the long-time average of a drag-related cost function, the energy of the velocity fluctuations in the wake of the cylinder. Before designing the control, the governing partial differential equations are projected on a finite-dimensional subspace obtained from Proper Orthogonal Decomposition,~\cite{Holmes:93}, to obtain a compact set of ordinary differential equations. Since the reduced system has a quadratic polynomial nonlinearity in the state variables, the proposed SOS-based control design can be applied. The performance of the proposed controller is assessed by closed-loop simulations of the reduced-order model as well as via direct numerical simulation (DNS).

The paper is organized as follows. Section~\ref{seq:Background} presents an introduction to SOS and its application to bound estimation of long-time average cost for uncontrolled systems. Bound optimization for controlled polynomial systems is considered in Section~\ref{seq:Control}. Section~\ref{seq:GeneralControlCase} extends the result obtained in Section~\ref{seq:Control} to more general scenarios. Then, feedback control design for the wake flow is addressed in Section~\ref{seq:Example}. Section~\ref{seq:Conclusion} concludes the work.

\section{Background}\label{seq:Background}

In this section, SOS of polynomials and a recently-proposed method of obtaining rigorous bounds of long-time average cost for uncontrolled polynomial systems are introduced.

\subsection{SOS of polynomials}

SOS techniques have been frequently used in the stability analysis and controller design for many kinds of systems, e.g., constrained ordinary differential equation systems~\cite{An:02}, hybrid systems~\cite{An:05}, time-delay systems~\cite{An:04}, and partial differential equation systems~\cite{Pa:06,Yu:08,GC:11}. These techniques help to overcome the common drawback of approaches based on Lyapunov functions: before~\cite{Pr:02}, there were no coherent and tractable computational methods for constructing Lyapunov functions.

A multivariate polynomial $f({\bf x}), ~\bfx\in {\mathbb R}^n$ is a SOS, if there exist polynomials $f_1({\bf x}), \cdots,$ $ f_m({\bf x})$ such that
\[
f({\bf x})=\sum_{i=1}^mf_i^2({\bf x}).
\]
If $f({\bf x})$ is a SOS then $f({\bf x})\ge 0, \forall{\bf x} $. In the general multivariate case, however,
$f({\bf x})\ge 0~  \forall \bf x$  does not necessarily imply that $f({\bf x})$ is a SOS. While being stricter, the condition that $f({\bf x})$ is SOS is much more computationally tractable than non-negativity~\cite{Par:00}. At the same time, practical
experience indicates that in many cases replacing non-negativity with the SOS property leads to satisfactory results.

In the present paper we will utilize the existence of efficient numerical methods and software~\cite{Pra:04,Lo:09} for solving the optimization problems of the following type~\cite{Pra:04}: minimize the linear objective function
\begin{equation}
{\bf w}^T{\bf c}
\label{linear}
\end{equation}
where $\bfw$ is the vector of weighting coefficients of the linear objective function, and ${\bf c}$ is a vector formed from the (unknown) coefficients of the polynomials $p_i({\bf x})$ for $i=1,2,\cdots, \hat{N}$ and SOS $p_i({\bf x})$ for $i=(\hat{N}+1),\cdots, {N}$, such that
\bea
\left\{
\begin{array}{c}
\displaystyle a_{0j}({\mathbf{x}})+\sum_{i=1}^Np_i({\mathbf{x}})a_{ij}({\mathbf{x}})=0, ~~~ j=1,2,\cdots, \hat{J}, \\
[1ex]
\displaystyle a_{0j}({\mathbf{x}})+\sum_{i=1}^Np_i({\mathbf{x}})a_{ij}({\mathbf{x}})  \mbox{~is SOS,~}~~ j=(\hat{J}+1),\cdots, {J}.~
\end{array}
\right.
\label{c6}
\eea
In (\ref{c6}), $a_{ij}(\bfx)$ are polynomials with given constant coefficients.

The lemma below provides a sufficient condition to test inclusions of sets defined by polynomials and is frequently used for feedback controller design in Section~\ref{seq:Control}. It is a particular case of Positivstellensatz~\cite{Po:99} and is a generalized $S$-procedure~\cite{Ta:06}.

\begin{lemma}
Consider two sets of $\bfx$,
\beas
{\mathcal S}_1&\eqdef& \left\{\mathbf{x}\in {\mathbb R}^n~|~h(\mathbf{x})=0, f_1(\mathbf{x})\ge 0, \cdots, f_r(\mathbf{x})\ge 0\right\}, \\
[1ex]
 {\mathcal S}_2&\eqdef& \left\{\mathbf{x}\in {\mathbb R}^n~|~f_0(\mathbf{x})\ge 0\right\},
\eeas
where $f_i(\mathbf{x}), i=0,\cdots, r$ and $h(\mathbf{x})$ are scalar polynomial functions. The set inclusion ${\mathcal S}_1\subseteq {\mathcal S}_2$ holds if there exist a polynomial function $S_0(\mathbf{x})$ and SOS polynomial functions $S_i(\mathbf{x}), i=1,\cdots, r$ such that
\beas
f_0(\mathbf{x})+S_0(\mathbf{x})h(\mathbf{x}) -\sum_{i=1}^rS_i(\mathbf{x})f_i(\mathbf{x})~~\mbox{~is a SOS.}
\eeas
\end{lemma}

\subsection{Bound estimation of long-time average cost for uncontrolled systems}\label{UncontrolledBound}

For the convenience of the reader we outline the method of obtaining bounds for long-time averages proposed in~\cite{Ph:14}.
Consider a dynamical system
\bea
\dot{\mathbf{x}}=\mathbf{f}(\mathbf{x}), ~\bfx\in {\mathbb R}^n,
\label{sys1}
\eea
where $\dot{\mathbf{x}}\eqdef d\mathbf{x}/dt$ and assume that the trajectories $\bfx(t)$ are uniformly bounded as $t\rightarrow \infty$ regardless of the initial condition $\bfx_0$.
Let $\Phi(\bfx)$ be the cost function. Suppose there exist a constant $C$  and a differentiable function $V(\bfx)$ such that
\bea
\dot{V}+\Phi-C\le 0.
\label{intro_1}
\eea
Due to the uniform boundedness of $\bfx(t)$, $V$ is also bounded as $t\rightarrow \infty$, so time averaging (\ref{intro_1}) gives
\[
\bar{\Phi}\le C,
\]
where $\bar{\Phi}$ refers to the long-time average of $\Phi$, namely,
\[
\bar{\Phi}=\lim_{T\rightarrow \infty}\frac{1}{T}\int_0^T\Phi({\mathbf{x}}(t))\,\mathrm{d}t.
\]
Hence, an upper bound of $\bar{\Phi}$ can be obtained by minimizing $C$ over $V$ under the constraint (\ref{intro_1}), or equivalently
\beas
\mathbf{f}\cdot\nabla_{\mathbf{x}} V+\Phi-C\le 0, ~\forall \bfx,
\eeas
where $\nabla_{\mathbf{x}} V$ denotes the gradient of $V$ with respect to $\bfx,$ so that $\dot{V}=\mathbf{f}\cdot\nabla_{\mathbf{x}} V.$
When the vector field $\bff$, the cost function $\Phi$, and the tunable function $V$ are restricted to polynomials, it suffices to solve the following SOS optimization problem:
\begin{equation}
\left.
\begin{array}{lc}
&\displaystyle\min_{V,~C}~C \\
[2ex]
&\mbox{s.t.} ~~-\left( \mathbf{f}\cdot\nabla_{\mathbf{x}} V+\Phi-C\right) \mbox{~is~SOS},
\end{array}
\right.
\label{SOS}
\end{equation}
which is a special case of (\ref{linear})-(\ref{c6}).

A better bound might be obtained by removing the requirement for $V$ to be a  polynomial and replacing (\ref{SOS}) with the requirement of non-negativeness. However, the resulting problem could be difficult, since the classical algebraic-geometry problem of verifying positive-definiteness of a general multivariate polynomial is NP-hard~\cite{An:02,An:05}. Notice that while  $V$ is similar to a Lyapunov function in stability analysis, it is not required to be positive-definite. Notice also that a lower bound of any long-time average cost of the system (\ref{sys1}) can be analyzed in a similar way, by reversing the sign of the inequality.

\begin{remark}\label{RemarkOnBoundedness}
For many systems the boundedness of the system state immediately follows from energy consideration. In general, if the system state is bounded this can often be proven using the SOS approach.
As an example, let $\mathcal{D}=\{\mathbf{x}~|~0.5\mathbf{x}^T\mathbf{x}\le \beta < \infty\}.$
Then $\mathcal{D}$ is a global attractor provided that 
\bea
\mathbf{x}^T\dot{\mathbf{x}}=\mathbf{x}^T\mathbf{f}(\mathbf{x})\le -(0.5\mathbf{x}^T\mathbf{x}-\beta) \quad \forall \mathbf{x}.
\label{SOS1}
\eea
Introducing a tunable polynomial $S(\mathbf{x})$ satisfying $S(\mathbf{x})\ge 0 ~\forall \mathbf{x}\in{\mathbb R}^n$, by Lemma~1, (\ref{SOS1}) can be relaxed to
\[
\left\{
\begin{array}{c}
-\left(\mathbf{x}^T\mathbf{f}(\mathbf{x})+S(\mathbf{x})(0.5\mathbf{x}^T\mathbf{x}-\beta)\right)\mbox{~is ~SOS},  \\
[1ex]
S(\mathbf{x})\mbox{~is ~SOS}.
\end{array}
\right.
\]
Minimization of upper bound of long-time average cost for systems that have unbounded global attractor is usually meaningless, since the cost itself could be infinitely large.
\end{remark}

\section{Expensive control of long-time averages of polynomial systems}\label{seq:Control}

In this section, the LTACC of polynomial systems is formulated first. 
Such systems may describe a wide variety of dynamics~\cite{Va:01} or approximate a system defined by an analytical vector field~\cite{Va:02}. A polynomial system can therefore yield a reliable model of a dynamical system globally or in larger regions than the linear approximation in the state-space~\cite{Va:03}. Finite-dimensional representations of incompressible fluid flows, as the example discussed in Section~\ref{seq:Example}, can be recast precisely in this form.

Then, the idea of the expensive controller design is presented to resolve the non-convexity in SOS optimization. 

\subsection{Problem formulation}\label{sec:problem formulation}

Consider the system
\bea
\dot{\mathbf{x}}=\mathbf{f}(\mathbf{x})+\mathbf{G}(\mathbf{x})\mathbf{u}
\label{sys}
\eea
where $\mathbf{f}(\mathbf{x}): {\mathbb R}^n\rightarrow {\mathbb R}^n$ and $\mathbf{G}(\mathbf{x}): {\mathbb R}^n\rightarrow {\mathbb R}^{n\times m}$ are polynomial functions of the system state $\bfx$. The control $\bfu\in {\mathbb R}^m$ is assumed to be a polynomial vector function of the system state $\bfx.$
 
 The cost function is a time average:
\[
\bar{\Phi}=\lim_{T\rightarrow \infty}\frac{1}{T}\int_0^T\Phi({\mathbf{x}}(t),\mathbf{u}({\mathbf{x}}(t)))\,\mathrm{d}t,
\]
where ${\mathbf{x}}$ is the closed-loop solution of the system (\ref{sys}) associated with the control $\mathbf{u}$, and 
the continuous function $\Phi$ is a given polynomial in ${\mathbf{x}}$ and $\mathbf{u}$.

Here, the existence of the upper bound of $\bar{\Phi}$ is assumed. This assumption holds true for many systems, for which the long-term behavior is determined by attractors. Note that even if the uncontrolled system has a global attractor, the closed-loop system can become unbounded, that is some of its trajectories can escape to infinity. In that case imposing additional constraints of the type 
(\ref{SOS1})  on the controlled system  might resolve the issue.

In~\cite{Ph:14} it was proposed to seek the control $\bfu$ minimizing upper bound $C$ of the cost $\bar{\Phi}.$ 
Similar to (\ref{SOS}), this reduces to the following SOS optimization problem:
\begin{equation}
\left.
\begin{array}{lc}
&\displaystyle \min_{\mathbf{u}, V, C}~~ C \\
[1ex]
&\mbox{s.t.} -\left((\mathbf{f}(\mathbf{x})+\mathbf{G}(\mathbf{x})\mathbf{u}(\mathbf{x}))\cdot\nabla_{\mathbf{x}} V(\mathbf{x})+\Phi(\mathbf{x}, \mathbf{u}(\mathbf{x}))-C\right) ~\mbox{is SOS}.
\end{array}
\right.
\label{eqn:optimization}
\end{equation}
Since the control input $\bfu$ and the decision function $V,$ both of which are tunable, enter~(\ref{eqn:optimization}) nonlinearly,~(\ref{eqn:optimization}) is not convex.


If the control is expensive, the optimal solution $\mathbf{u}$ should be small. This suggests a Taylor expansion of $\Phi(\mathbf{x}, \mathbf{u})$ with respect to $\mathbf{u}.$ In many cases it is natural to assume that  any non-zero control is more expensive than no control, that is that that $\Phi(\mathbf{x}, \mathbf{u}(\mathbf{x}))$ has a minimum at $\mathbf{u}=0$ for  any fixed $\mathbf{x}.$ In this case the Taylor expansion has a zero linear term. Neglecting the cubic and the higher-order terms then leads to considering the cost function of the form

\bea
\Phi(\mathbf{x}, \mathbf{u})=\Phi_0(\mathbf{x})+\dfrac{\mathbf{u}^T \mymtrx{\Psi}(\bfx)\mathbf{u}}{\epsilon},
\label{new_cost}
\eea
where $\mymtrx{\Psi}$ is a symmetric positive-definite matrix of polynomials of $\bfx$, and $\epsilon$ is a small parameter, which makes the resulting problem fall in the class of expensive control~\cite{Chen:07}.
%
%

We are seeking to find the main terms of the asymptotic expansions of $C$ and $\mathbf{u}_{\text{opt}}(\bfx)$ as $\epsilon\to0.$

%
%

\subsection{Design of small-feedback controller}

We look for a controller in the form
\bea
\mathbf{u}(\mathbf{x},\epsilon)=\epsilon\mathbf{u}_1(\bfx)+O(\epsilon^2).
\label{cc}
\eea
We expand $V$  and $C$ in $\epsilon$:
\bea
V(\mathbf{x},\epsilon)&=&V_0(\bfx)+\epsilon V_1(\bfx)+O(\epsilon^2), \label{LF0} \\
C(\epsilon)&=&C_0+\epsilon C_1+O(\epsilon^2).
\label{LF1}
\eea
 Define
\bea
F(V,\bfu,C)\eqdef (\mathbf{f}(\mathbf{x})+\mathbf{G}(\mathbf{x})\mathbf{u})\cdot \nabla_{\mathbf{x}} V+\Phi(\mathbf{x},\mathbf{u})-C,
\label{cc1}
\eea
so that the constraint in the optimization problem~(\ref{eqn:optimization}) is  $F(V,\bfu,C)\le0\ \forall \mathbf{x}.$  
Substituting (\ref{cc}), (\ref{LF0}), and (\ref{LF1}) into (\ref{cc1}), we have
\beas
\begin{array}{ll}
F(V,\bfu,C)=&F_0(V_0,C_0)+\epsilon F_1(V_0,V_1,\bfu_1,C_1)+O(\epsilon^2),
\end{array}
\eeas
where
\beas
\begin{array}{ll}
F_0=&\bff\cdot \nabla_{\bfx} V_0+\Phi_0-C_0, \\
F_1=&\bff\cdot \nabla_{\bfx} V_1+\mathbf{G} \bfu_1\cdot \nabla_{\bfx} V_0+\mathbf{u}_1^T \mymtrx{\Psi}(\bfx)\mathbf{u}_1-C_1.
\end{array}
\eeas

Substituting these expansions into~(\ref{eqn:optimization}) and taking the limit as $\epsilon\to0$ leads to the problem of finding the best bound for the uncontrolled case: 
\beas
O_0: ~~\min_{V_0, C_0} C_0, \mbox{~~s.t.~~} -F_0(V_0,C_0)\mbox{~~is SOS},
\eeas 
which we will call $O_0$ problem.

Denote the optimal $C_0$ by $C_{0,SOS}$ and the associated $V_0$ by $V_{0,SOS}$.

Let $V_0=V_{0,SOS}$ in $F_1$, and then consider the following optimization problem:
\beas
&O_1: \begin{array}{c}
\displaystyle \min_{V_1,\bfu_1,S_0, C_1} C_1, \\
[2ex]
\mbox{~~s.t.~~}
\begin{array}{c}
-F_1(V_{0,SOS},V_1,\bfu_1,C_1) \\
+S_0(\bfx)F_0(V_{0,SOS},C_{0,SOS}) \mbox{~is SOS},
\end{array}
\end{array}
\eeas
where $S_0$ is a tunable polynomial function of $\bfx$ of fixed degree corresponding to the $S$-procedure.
Here the non-negativity requirement of $S_0$ is not imposed. By Lemma~1, this can be understood as that the non-negativity constraint $-F_1(V_{0,SOS},V_1,\bfu_1,C_1)\ge 0$ is imposed only for $\bfx$ such that $F_0(V_{0,SOS},C_{0,SOS})=0$.
Denote the optimal $C_1$ by $C_{1,SOS}$ and the associated $V_1$ and $\bfu_1$ by $V_{1,SOS}$ and $\bfu_{1,SOS}$, respectively.

Note that since $F_1$ includes a quadratic term in $\bfu_1$, the optimisation problem $O_1$ is not convex.
This can be overcome using the idea of the Schur complement. In the current setting it amounts to introducing an additional variable $\mathbf{z}$ and replacing  $F_1$ with $F'_1=F_1- (\mathbf{u}_1-\mathbf{z})^T \mymtrx{\Psi}(\bfx)(\mathbf{u}_1-\mathbf{z}).$ The quadratic term then cancels out. Since $\mymtrx{\Psi}$ is positive-definite, the constraint  $-F_1+S_0F_0>0$ $\forall \bfx$ is equivalent to the constraint  $-F'_1+S_0F_0>0$ $\forall \bfx, \mathbf{z}.$ This reduces $O_1$ to an optimisation problem of the type (\ref{linear}-\ref{c6}).

The solution of $O_1$ always gives $C_1\le0$, since $C_1=0$ for $\bfu_1=0.$ If the optimal $C_1<0$ then, as it follows from (\ref{LF1}), the solution gives a controller reducing the bound of the time-averaged cost in an asymptotic, as $\epsilon\to0,$ sense. The following theorem shows that in fact the controller $\epsilon \bfu_1$ reduces the bound for a sufficiently small but finite $\epsilon$ too.

\begin{theorem}
Assume that $C_{1,SOS}<0$, obtained by solving $O_1$. Then, for any $\kappa\in(0,1)$, there exists $\epsilon_{\kappa}>0$ such that
applying the following first-order small-feedback controller to the system (\ref{sys}),
\bea
\bfu_{SOS}=\epsilon \bfu_{1,SOS}, ~~0<\epsilon<\epsilon_{\kappa},
\label{new3}
\eea
will yield an upper bound of the long-time average cost $\bar{\Phi}$,
\beas
C_{\kappa,SOS}\eqdef C_{0,SOS}+\epsilon\kappa C_{1,SOS}.
\eeas
Clearly, $C_{\kappa,SOS}<C_{0,SOS}$.
\end{theorem}

{\it Proof}. See {\bf Appendix A}.

\commentout{
\begin{remark}
It is worth noticing that the linear approximation on the controller does not mean that our controller design and analysis are conducted in a non-rigorous way. The truncated small-feedback controller would be effective if it leads to a better (lower) bound of the long-time average cost.
\end{remark}

\begin{remark}
By specifying the structure of controller as in (\ref{cc}), the non-convexity in solving the optimization problem (\ref{objective})-(\ref{optimization}) has been avoided by solving $O_0, O_1, \cdots$ in sequence.
During the process, all the involved decision variables are optimized sequentially, but not iteratively as in other methods~\cite{Zh:07,Zh:09,Ng:11}.
Nevertheless, due to the existence of quadratic terms of tuning variables, e.g., $\left|\bfu_1\right|^2$ in $F_1$, linear SDP solvers are still not applicable. In Section~\ref{seq:Example}, we will adopt the Schur complement formula to overcome this issue.
\end{remark}
}

\begin{remark}

Once the controller has been specified as in (\ref{new3}), with $\epsilon$ and $\bfu_{1,SOS}$ given, the upper bound $C$ and the corresponding $V$ can be obtained by solving a smaller optimization problem 
\beas
O_{\epsilon}: \begin{array}{c}
\displaystyle  \min_{V, \epsilon, C} C, ~~~ \\
 [1ex]
\mbox{s.t.}~~-F(V,\epsilon \bfu_{1,SOS},C)\mbox{~is~ SOS}.
\end{array}
\eeas
 The smaller size of the problem will allow searching over polynomial $V$ of higher order. The  problem can be further relaxed by utilizing the knowledge that $-F_0$ and $-F_1+S_0F_0$ are sums of squares. This might allow getting a better bound.
The effect of the value of $\epsilon$  on the upper bound of $\bar{\Phi}$ can be investigated by trial and error. We will follow this route in Section~\ref{seq:Example}.

\end{remark}

\section{More general case}\label{seq:GeneralControlCase}
\label{generalcase}

In practice, some dynamical systems (as the one discussed in Section~\ref{seq:Example}) might have more complicated dynamics described by the following form
\bea
\dot{\bfx}=\bff(\bfx)+\bfG(\bfx)\bfu+\bfH(\bfx)\dot{\bfu},
\label{gensys}
\eea
where the time rate of change of the control input $\dot{\bfu}$ also enters the dynamics via the polynomial gain function $\bfH(\bfx): {\mathbb R}^n\rightarrow {\mathbb R}^{n\times m}$.

With a minor modification, the small-feedback controller design scheme is still applicable for (\ref{gensys}). Similar to (\ref{cc1}), define
\bea
F(V,\bfu,C)\eqdef (\bff(\bfx)+\bfG(\bfx)\bfu+\bfH(\bfx)\dot{\bfu})\cdot \nabla_{\bfx} V+\Phi(\bfx,\bfu)-C.
\label{gencc1}
\eea
Noticing the structure of the controller (\ref{cc}), it is easy to see that
\bea
\dot{\bfu}&=&\frac{\mathrm{d}\bfu}{\mathrm{d}\bfx} \dot{\bfx}=\frac{\mathrm{d}\bfu}{\mathrm{d}\bfx}\left(\bff(\bfx)+\bfG(\bfx)\bfu+\bfH(\bfx)\dot{\bfu} \right) \nonumber \\
&=&\epsilon \frac{\mathrm{d}\bfu_1}{\mathrm{d}\bfx} \bff(\bfx) +O(\epsilon^2).
\label{u_derivative}
\eea
Then, substituting (\ref{cc}), (\ref{LF0}), (\ref{LF1}), and (\ref{u_derivative}) into (\ref{gencc1}) renders to
\beas
F=F_0(V_0,C_0)+\epsilon F_1(V_0,V_1,\bfu_1,C_1)+O(\epsilon^2),
\eeas
where
\beas
F_0&=&\bff\cdot \nabla_{\bfx}V_0+\Phi_0-C_0,
\\
F_1&=&\bff\cdot \nabla_{\bfx}V_1+( \bfG\bfu_1+\bfH\frac{\mathrm{d}\bfu_1}{\mathrm{d}\bfx}\bff )\cdot\nabla_{\bfx}V_0+\mathbf{u}_1^T \mymtrx{\Psi}(\bfx)\mathbf{u}_1-C_1.
\eeas
The corresponding SDP problems may be solved via SOS optimization to obtain the controller that minimizes the bound of the long-time average cost of the controlled system (\ref{gensys}).

\section{Illustrative example}\label{seq:Example}

As an illustrative example of the methodology proposed in this paper, we consider the problem of synthesizing a state-feedback controller to manipulate the motion of a fluid flow. Fluids are a prominent example where stabilization of the laminar flow, (usually associated with low drag, low unsteadiness, low aerodynamic-induced vibrations, etc.), might not be possible in practice, extremely difficult to achieve or not worth the additional complexity and cost of the required flow sensing/actuation system. In addition, as fluid flows are often turbulent, exhibiting chaotic fluctuations over a disparate range of time and spatial scales, long-time averages of key engineering quantities are often of interest. Hence, nonlinear design techniques aiming at reducing, even by a small amount, such averages are attractive and relevant. For instance, a few percent reduction of the mean drag of an aircraft, associated with the control of the turbulent flow over its surface, would entail significant economic benefits and would have a large impact on its environmental footprint~\cite{kim:07}.

The problem of mitigating vortex shedding in the two-dimensional flow past a circular cylinder via a controlled rotary motion of the cylinder is often adopted as a benchmark to test strategies and algorithms of of flow control, as it is a relevant paradigm of separated flows past bluff bodies. The full stabilization of the wake past the cylinder, i.e. the complete suppression of vortex shedding, has been achieved only at very low Reynolds numbers. Linear design methods have been also considered extensively for this problem, see for example~\cite{Ro:14, carini:15, illingworth:14} and references therein, to stabilize the unstable steady solution of the equation. However, these methods seem to become ineffective when the controller needs to control the fully-developed nonlinear regime. Nonlinear optimal control, in the receding-horizon setting, is considered in~\cite{Pr:02} and more recently in~\cite{flinois:15}. In this setting, the control in the form od the angular velocity of the cylinder over a finite-horizon is found in real time from an expensive optimisation procedure, involving repeated solutions of the governing partial differential equations and of their adjoints,~\cite{abergel:90, Bewley:01}. In~\cite{flinois:15}, stabilization of the flow was achieved at higher Reynolds numbers than in previous works, although the method is computationally expensive and the controllability of the system will worsen at larger Reynolds number.

In this section  we will test the method proposed in the previous section on this benchmark problem.

\commentout{
\begin{figure}[htbp]
	\centering
	\includegraphics[width=0.5\textwidth]{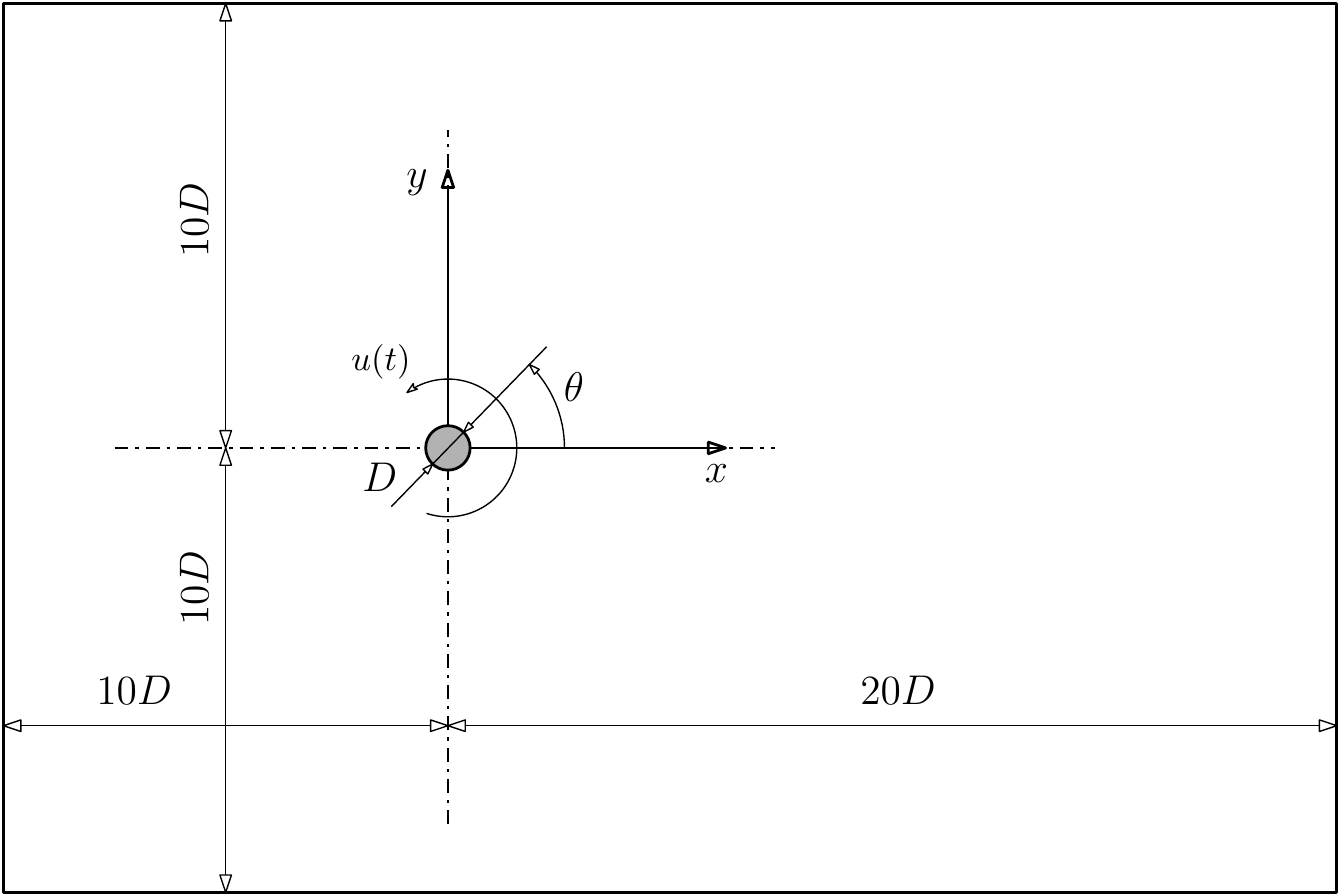}
	\caption{Sketch of the computational domain used for direct numerical simulation of control of two-dimensional vortex shedding past a circular cylinder.}
	\label{fig:shedding}
\end{figure}
}

\subsection{Problem formulation}
The flow of the viscous incompressible fluid past the cylinder is described by the Navier-Stokes and continuity equations
\begin{equation}\label{eq:ns}
	\frac{\partial \mathbf{v}}{\partial t} +\mathbf{v} \cdot \nabla \mathbf{v} 
	= - \nabla p + \frac{1}{Re}
	\nabla^2
	 \mathbf{v}, 
	\qquad
	\nabla \cdot \mathbf{v} 
	= 0,
	\end{equation}
where $\mathbf{v}(\mathbf{x}, t)$ is the velocity vector field defined on the two-dimensional Cartesian space $\mathbf{x}=(x, y)$, and $p$ is the pressure. The Reynolds number is $Re=u_\infty D /\nu$, where $D$ is the cylinder diameter, $u_\infty$ is the free stream velocity and $\nu$ is the kinematic viscosity of the fluid. Velocities, pressure, lengths and times have been made non-dimensional using $u_\infty$ and $D$, or combinations thereof. All numerical results of this section were obtained for $Re=100.$ 

\commentout{
 In equation (\ref{eq:ns}), $\nabla = (\partial/\partial x, \partial/\partial y)$ is the nabla operator, whereas $\Delta = \partial^2/\partial x^2 + \partial^2/\partial y^2$ is the Laplacian operator.
}

Flow actuation is performed via the boundary, via controlled rotary motions of the cylinder. The scalar time-dependent control input $u(t) \in \mathbb{R}$ is the tangential velocity normalised with the free stream speed. To close the feedback loop, we assume that full information on the velocity vector field is available, as the focus of this paper is on the control design method.

Direct numerical simulation (DNS) of the flow is performed using the OpenFoam package~\cite{of:98}, on an unstructured triangular mesh. The computational domain $\Omega$ is the rectangle extending 10 and 20 diameters upstream and downstream the cylinder, respectively, and 20 diameters wide. Free-slip boundary conditions are imposed on the upper and lower horizontal boundaries. On the inlet, the Dirichlet condition $\mathbf{v} = (1, 0)$ is imposed, whereas the Neumann condition $\partial p/\partial x = 0$ is used for pressure. At the outlet, good numerical results have been obtained by using the condition $\partial \mathbf{v}/\partial x = (0, 0)$, whereas the Dirichlet condition $p=0$ was set to fix uniquely the pressure field. On the cylinder surface, the velocity component normal to the cylinder is set to zero,  the tangential velocity is specified by the time-dependent control input $u(t)$, and a standard zero normal pressure gradient condition is specified for the pressure.  The nondimensional time step was set to $\Delta t = 0.005$, to limit the maximum value of the CFL number below one. Validation and grid convergence studies have been performed but are not reported in the present paper as the numerical method and the flow problem are rather standard.

In what follows, the inner product between two vector fields $\mathbf{v}$ and $\mathbf{w}$, defined as
\[
	\langle \mathbf{v}, \mathbf{w}\rangle = \int_\Omega \mathbf{v} \cdot \mathbf{w} ~\mathrm{d}\Omega,
\]
will be used. Such a definition implies that the norm of a vector field is $\|\mathbf{v}\| = \langle\mathbf{v}, \mathbf{v}\rangle^{1/2}$.

\subsection{Low-order modelling}
A finite-dimensional description of the dynamics, in the form of a set of first-order ordinary differential equations with right-hand side that is a polynomial function in the state variables and in the control, as in equation (\ref{gensys}), is required for control design. However, the system (\ref{eq:ns}) is a partial differential equation or, after discretisation, a high-dimensional system that will also be referred to as the full-order system in the following. In this paper, the Proper Orthogonal Decomposition and Galerkin projection,~\cite{Holmes:93}, are adopted to obtain a reduced-order description, a reduced order model (ROM) that compactly describes the actuated dynamics of the wake and that can be used for control design. In what follows, we report for the sake of completeness the modelling strategy adopted in this paper, which follows the works~\cite{Gr:99, be:05}. The interested reader is referred to these works and references therein for a more detailed description of this strategy.

A truncated Galerkin expansion of the velocity field defined by the ansatz
\begin{equation}\label{eq:ansatz}
	\mathbf{v}^N(\mathbf{x}, t) = \overline{\mathbf{v}}(\mathbf{x}) + u(t)\mathbf{v}_c(\mathbf{x}) + \sum_{i=1}^N a_i(t)\mathbf{v}_i(\mathbf{x})
\end{equation}
is first introduced. Here, the velocity field $\mathbf{v}^N$ is decomposed into a sum of a mean flow $\overline{\mathbf{v}}$ satisfying homogeneous boundary conditions on the cylinder, a ``control function'' $u(t)\mathbf{v}_c(\mathbf{x})$, (see e.g.~\cite{Gr:99,kasna:08}) used to lift the time-dependent inhomogeneous boundary conditions on the oscillating cylinder surface and to include control via the boundary in the dynamic model, and a weighted sum of $N$ basis functions $\mathbf{v}_i(\mathbf{x})$ forming an orthonormal set. A radially-symmetric control function $\mathbf{v}_c(\mathbf{x})$, with circumferential velocity decaying as $e^{-8(r-0.5)}$, was employed.

The snapshot variant of POD,~\cite{sirovich:87}, is used to derive the basis functions $\mathbf{v}_i(\mathbf{x})$. Following the procedure described in~\cite{be:05}, the first set of snapshots of the velocity vector field, $\mathcal{V} = \{\mathbf{v}(\mathbf{x}, t_k)\}_{k=1}^M$, is sampled from a direct numerical simulation in which the angular motion of the cylinder is driven by a random actuation signal, with the idea of exciting transient flow structures and obtaining a richer snapshot set. The signal is obtained from samples of a zero-mean Gaussian distribution, and it is then filtered such that its power spectrum has zero energy outside the band of reduced frequency $fD/u_\infty = [0.1, 0.25]$, by setting to zero the appropriate entries of its Fourier transform, before transforming back to the real space. The amplitude of the filtered signal is then modulated by a mode with reduced frequency $fD/u_\infty = 0.005$, in order to actuate the flow at different intensities, and it is then normalised to have unitary maximum magnitude, resulting in a standard deviation equal to about 0.25. The total duration of this simulation is $T=1000$, about 150 oscillation cycles of the uncontrolled flow, and a total of $M=900$ snapshots is sampled, from $t \ge 100$, at intervals of 1 non-dimensional time unit.

The time-dependent, inhomogeneous boundary conditions on the cylinder are then removed from the snapshots by subtracting, with appropriate amplitude, the control function, obtaining the set
\[
	\displaystyle \mathcal{V}' = \{\mathbf{v}^h(\mathbf{x}, t_k) = \mathbf{v}(\mathbf{x}, t_k) - u(t_k)\mathbf{v}_c(\mathbf{x})\}_{k=1}^M.
\]
The arithmetic average is then computed
\[
	\overline{\mathbf{v}}(\mathbf{x}) = \frac{1}{M}\sum_{k=1}^M \mathbf{v}^h(\mathbf{x}, t_k)
\]
and it is used as the mean flow for the ansatz (\ref{eq:ansatz}). Finally, the snapshot set
\[
	\mathcal{V}'' = \{ \mathbf{v}^h(\mathbf{x}, t_k) - \overline{\mathbf{v}}(\mathbf{x})\}_{k=1}^M
\]
is used for the POD algorithm, yielding the basis functions $\mathbf{v}_i(\mathbf{x}),\; i=1, \ldots, M$.

We selected the first $N=9$ POD modes for the Galerkin projection, as a compromise between the accuracy of the model to resolve the actuated dynamics of the wake and the computational costs of the solution of the SOS problems described in the previous sections. Furthermore, a shift mode is added to improve transient dynamics over changes in the base flow,~\cite{noack:03}. The shift mode is a basis function spanning the direction from the mean flow $\overline{\mathbf{v}}(\mathbf{x})$ to the unstable, steady and symmetric solution $\mathbf{v}_0(\mathbf{x})$ of (\ref{eq:ns}), obtained numerically as the steady solution on the upper half domain, with free-slip boundary condition on the symmetry plane. The shift mode is constructed as
\[
	\mathbf{v}_\Delta(\mathbf{x}) = \frac{\overline{\mathbf{v}} - \mathbf{v}_0}{\|\overline{\mathbf{v}}-\mathbf{v}_0\|}
\]
and it is made orthogonal to the remaining nine POD modes using a Gram-Schmidt procedure.

Galerkin projection is performed by inserting the expansion (\ref{eq:ansatz}) in (\ref{eq:ns}), and setting the inner product with each of the modes in turn to zero. Neglecting the small contribution arising from the projection onto the pressure gradient field, as commonly done for this fluid flow (e.g.~\cite{be:05, noack:03}), results in the nonlinear reduced-order model (ROM):
\begin{equation}\label{eq:ode-sys}
\begin{array}{rl}
\displaystyle \dot{a}_i = &c_i + {\sum_{j=1}^NL_{ij} a_j}  + {\sum_{j=1}^N \sum_{k=j}^N N_{ijk}a_j a_k} \\
&+ {m_i   \dot{u}  + e_i u + b_i u^2 + \sum_{j=1}^N F_{ij} a_j u}, \quad i=1, \ldots N.
\end{array}
\end{equation}
The definitions of the coefficients $c_i, L_{ij}, N_{ijk}, m_i, e_i, b_i, F_{ij}$ arising from the projection are given in Appendix B. Because of the particular choice of the function $\mathbf{v}_c$, the coefficients of the quadratic term $b_i$ vanish identically. The matrices  $L_{ij},$ $ N_{ijk},$ and $ F_{ij}$ are dense, which makes it difficult to include them in the present text.
In vector form, the system (\ref{eq:ode-sys}) is
\[
\dot{\mathbf{a}}=\mathbf{c}+L\mathbf{a}+\mathbf{N}(\mathbf{a})\mathbf{a}+\mathbf{e}u+F\mathbf{a}u+\mathbf{m}\dot{u}
\]
which is in the form given by equation (\ref{gensys}), discussed in Section~\ref{generalcase}. We assume that the full information on the system state is available. In the direct numerical simulation, the state vector is obtained from projection of the POD basis functions on the solution, i.e.
\[
a_i(t) = \langle \mathbf{v}_i(\mathbf{x}), \mathbf{v}(\mathbf{x}, t) - \overline{\mathbf{v}}(\mathbf{x})\rangle.
\]
In a physical experiment, the system state could be estimated by designing a suitable state observer, based on wall or field measurements~\cite{Ro:14}. The ROM (\ref{eq:ode-sys}) is integrated numerically using a standard fourth-order Runge-Kutta method, with the time step $\Delta t = 0.001$.

For this flow problem, we select the term $\Phi_0$ in (\ref{new_cost}) to be the energy of the system
\[
	\Phi_0(\mathbf{a}) = \frac{1}{2} \mathbf{a}^T\mathbf{a},
\]
similarly to other works, e.g.~\cite{be:05}. Physically, $\Phi_0$ represents the domain integral of the kinetic energy of the velocity fluctuations resolved by the ansatz (\ref{eq:ansatz}). As a result, reduction of this quantity will result in a reduction of the wake unsteadiness associated with vortex shedding.

The ten-mode dynamical system obtained from projection does not represent very accurately the dynamics of the full-order system, that is the numerical solution of the discretized governing equations (\ref{eq:ns}). In particular, the long-term behaviour, the stable limit cycle associated with vortex shedding, is not correctly represented. In fact, numerical integration of the ROM with $u(t) = 0$ shows that trajectories converge to a stable limit cycle with a long-time-averaged cost $\Phi_0$ about three times higher than that obtained from the projections of the basis functions on the long-term solution of the full-order system.
Because the long-time average cost of dynamical systems is usually based on the structure and type of the invariant sets of the system,  it is desirable to have a ROM whose long-term behaviour is as similar as possible to that of the full-order system. To this end, we apply a model calibration scheme, which has become the standard practice to correct the neglected effects of truncated modes on the resolved modes~\cite{couplet:07,sirisup:04}.

Following previous work (see e.g.~\cite{cordier:10}), a linear calibration term $L_{ij}^c$, with non-zero elements on the main, first upper and first lower diagonals is added to the linear term $L_{ij}$ in (\ref{eq:ode-sys}). This calibration term is obtained as the solution of an optimization problem in which the integral of the norm of the error between the ROM trajectory and the projection of the trajectory of the full-order system onto the selected subspace, obtained from a numerical simulation of the uncontrolled flow, is minimised. In more formal notation we solve numerically
\[
\displaystyle \underset{L_{ij}^c}{\text{min}} \displaystyle \int_{t_0}^{t_1} \| \mathbf{a}(t; t_0, L^c_{ij}) - \tilde{\mathbf{a}}(t; t_0) \|^2 \mathrm{d}t
\]
where $\tilde{a}_i$ are the projections of the DNS solution onto the POD basis functions, the notation $\mathbf{a}(t; t_0; L^c_{ij})$ denotes that the calibrated ROM is integrated in time from an initial condition $\mathbf{a}(t_0) = \tilde{\mathbf{a}}(t_0)$, with $u(t)=0$. The final time $t_1$ is such that $t_1-t_0$ amounts to 10 vortex shedding periods.

The calibrated ROM possesses a stable limit cycle closer to that of the full-order system. However, poor controllability of this reduced-order model was observed, as opposed to larger models that did not present this behaviour, suggesting that the rotary actuation of the cylinder affects via viscosity the large scale motion, i.e. the resolved modes, through linear/nonlinear interaction of the truncated modes. To mitigate this poor controllability, we also calibrate part of the coefficients associated with control in (\ref{eq:ode-sys}), in particular those associated with the term $e_i$ and $m_i$, against the numerical simulation used to obtain the POD modes, using a similar procedure as for the term $L_{ij}^c$. This resulted in the ROM suitable for control design for this particular flow.

\commentout{

\begin{figure*}[ht]
\begin{align*}
\mathbf{c}=&
\left(
\begin{array}{cccccccccc}
-0.0009&
-0.0007&
0.0026&
-0.0008&
0.0027&
0.0217&
0.0729&
-0.0145&
0.1548&
-0.0675
\end{array}\right)^T,\\
\mathbf{m}=&
\left(
\begin{array}{cccccccccc}
    1.4137&
    0.0758&
    0.1179&
   -0.5210&
   -0.5549&
    0.1610&
   -0.2358&
   -0.0809&
   -0.1262&
   -0.0092
\end{array}
\right)^T, \\
\mathbf{e}=&
\left(
\begin{array}{cccccccccc}
   -0.2927&
   -1.5626&
    0.4549&
    0.1830&
    0.8453&
    0.3336&
   -0.0164&
   -0.1801&
    0.1855&
   -0.0026
\end{array}
\right)^T,\\
L=&
\left(
\begin{array}{cccccccccc}
    0.1154 &  -1.0620 &   0.1175 &  -0.0760  &  0.0229 &   0.0495 &  -0.0178 &  -0.0010 &  -0.0074  & -0.0013\\
    1.0831 &  -0.1002 &  -0.0882  & -0.0912  & -0.0827 &   0.0208  & -0.0041 &   0.0057 &   0.0104 &   0.0006\\
   -0.0159 &   0.0148 &   0.0096   & 1.1302  &  0.2537 &   0.0415 &  -0.0325 &  -0.0108 &   0.0066  &  0.0008\\
    0.0096 &   0.0114 &  -1.0119  & -0.0303&   -0.1041 &  -0.1725 &   0.0725 &  -0.0058 &  -0.0355 &  -0.0028\\
    0.0114 &   0.0202 &  -0.1058 &   0.0359 &  -0.0140 &  -1.1971 &   0.2921 &  -0.0312 &  -0.0740 &  -0.0016\\
   -0.0119 &   0.0142 &  -0.0362 &   0.0654 &   1.0306 &  -0.2513 &  -0.1487 &   0.0245 &   0.0306 &   0.0184\\
    0.0027 &  -0.0062 &   0.0219 &  -0.0265 &  -0.2784 &   0.0863 &  -0.0296 &   0.5060  &  0.1038 &   0.0554\\
   -0.0020 &  -0.0003 &  -0.0033 &   0.0034 &   0.0215 &  -0.0262 &  -0.4299 &  -0.4085 &   0.1112 &  -0.0052\\
   -0.0024 &   0.0012 &  -0.0101 &   0.0127 &   0.0307 &  -0.0334 &  -0.1029  & -0.4871   &-0.4404 &   0.1325\\
    0.0000 &   0.0001 &  -0.0004 &   0.0009 &  -0.0000 &  -0.0111 &  -0.0359 &   0.0096  &  0.0157 &  -0.0462
\end{array}
\right),\\
F=
10^{-4}&\left(
\begin{array}{cccccccccc}
   -0.0292 &  -0.0400  & -0.0107  &  0.1832 &  -0.2406 &  -1.4190 &  -3.5245 &  -3.3165 &  -9.0367 &   1.3509\\
    0.1142 &  -0.0416&    0.1588 &  -0.1028 &   0.2272 &  -0.6118 &  -2.5986 &   1.4623 &  -3.3711  &  3.8802\\
   -0.1300 &   0.0075&   -0.1581 &   0.2442 &  -0.4277 &  -0.6356 &  -0.6264 &  -4.0601 &  -4.5169 &  -2.3645\\
   -0.1895 &   0.1452 &  -0.2918 &  -0.0322 &  -0.1752 &   1.7752 &   5.7059 &   1.8130 &  11.5645 &  -4.7837\\
    0.0760 &  -0.2019  &  0.2039  &  0.2465 &   0.0553  & -0.1981 &   0.4207 &  -5.0920 &  -3.4160  & -4.5864\\
    2.1898 &  -1.0912 &   2.8759 &  -0.9954 &   3.8102 &  -0.3282 &  -9.3164  & -3.4984 & -17.5911 &   6.8237\\
    6.0873 &  -2.7029 &   7.7632 &  -3.5162 &  11.4207 &   8.5054  &  0.0529  &  0.0287 &   6.4557 & -3.3207\\
    3.1637 &  -2.2968 &   4.7256 &  -1.0082 &   5.9223 &   2.3947  & -1.5398 &  -0.0801 &  -1.4423  &  0.9145\\
   13.5382 &  -6.8487 &  17.9216 &  -6.9518 &  25.3207  & 14.4743 &  -9.7152 &  -1.9404 &  -6.5295 &   2.9983\\
   -3.9835  &  0.6411 &  -4.1234 &   3.3042 &  -6.3723 &  -6.7127 &   3.1720 &  0.7003  & -0.5129  &  0.3353
\end{array}
\right).
\end{align*}
\noindent\rule{1\textwidth}{0.1pt}
\end{figure*}

}

\subsection{Calculation of the upper bound}
In what follows, the maximum degrees of polynomials $V$, $S$ and $u$ are denoted by $d_{V}$, $d_{S}$ and $d_{u}$ respectively. The SDP problem $O_0$ is solved first, to provide an estimate of the upper bound of the uncontrolled system. The least upper bound found by straightforward application of the method is $C_{0, SOS} = 167$, for a degree $d_{V_0}$ of the polynomial function $V$ equal to $4$. This value is considerably higher than the long-time averaged cost $\overline{\Phi}_0=3.07$  obtained from long-time numerical integration of the ROM (\ref{eq:ode-sys}) and by discarding initial transients before the trajectory has converged to the stable limit cycle associated with vortex shedding.

The likely reason for this discrepancy is the existence of a spurious invariant set, (a fixed point, a periodic orbit or another set), far away from the origin in the ten-dimensional phase space of the ROM, which, by design, can be expected to  approximate the full Navier-Stokes system only in the vicinity of the actual attractor. Such a spurious set, even unstable, would  affect negatively the calculation of the bound. 
\commentout{
Note that the stability of this set does not affect the calculations of the bound, and additional modification of the inequalities (\ref{optimization}) needs to be added to get rid of unstable set, as hinted in~\cite{Ph:14}. Application of the controller designed to reduce this large upper bound might result in an increase of the long-time averaged cost in numerical simulation.
}
This difficulty can be resolved  by confining the volume of the phase space where the polynomial inequalities involving the bounds need to be satisfied. For this, we consider the set described by the ball $\mathcal{B} = \{ \mathbf{a} \in \mathbb{R}^{10} : 1/2\bfa^T\bfa\leqslant 4\}$ that contains the periodic orbit associated with vortex shedding. To satisfy the polynomial inequality associated with the estimation of the bound in this ball rather than in the entire space  the $S$-procedure is used: a polynomial function in the state variables $S_0(\mathbf{a})$, of degree $d_{S_0}$ and with unknown coefficients as decision variables is introduced to modify the problem $O_0$ to
\beas
O_0': \begin{array}{c}
\displaystyle \min_{V_0,C_0}~~ C_0, \\
[1ex]
{~~s.t.~~} \left\{
\begin{array}{c}
-F_0(V_0,C_0)
-S_0(\bfa)(8-\bfa^T\bfa) \mbox{~~is SOS}, \\
[1ex]
S_0(\bfa)\mbox{~~is~ SOS}.
\end{array}
\right.
\end{array}
\eeas
For the case $d_{V_0}=4$ and $d_{S_0}=4$, the minimal upper bound we achieve is $C_{0,SOS}=3.07$, which is close to the value obtained from long-time integration of the ROM.

\subsection{Design of small-feedback controller}

In many situations the dependence of the cost function on the control law  is prescribed. 
However, in certain aeronautical applications the cost of implementing the control, that is the cost of sensors, actuators, and the associated overheads, such as for example the drag penalty due to the increased mass of the aircraft, while significant,  is independent of the control law.  At the same time, once the control device is installed, 
its operational cost, that is for example the additional fuel consumption needed to operate the controller, might be small as compared to the potential gain. In this situation the actual cost is in fact independent of the control law.  On the other hand, in practice the implementation cost would depend on the required magnitude of control, so that it is desirable to keep it small. With this situation in mind, an artificial penalty on the control law, dependent on a parameter, can be added to the cost function. The control law can then be designed, and the parameter can be varied afterwards to evaluate the usefulness of control. With this in mind we will select the particular form of the penalty term to be simple. Namely, we will assume that the cost function is
$\Phi=\Phi_0+||\bfu||^2/\epsilon,$ and where $\epsilon$ is an artificial expensive-control parameter, assumed to be small. 
A small-feedback controller $u=\epsilon u_1$ is to be designed to reduce the time-averaged cost $\bar{\Phi}.$ 

Similarly to $O_0'$, an additional $S$-procedure is applied to the problem $O_1$ to enforce the polynomial inequalities only in the ball $\mathcal{B}$. This leads to the modified problem
\beas
&O_1': \begin{array}{c}
\displaystyle \min_{V_1, u_1,S_0,S_1, C_1} C_1, \\
[1ex]
\mbox{~~s.t.~~}
\left\lbrace
\begin{array}{c}
-F_1(V_{0,SOS},V_1,u_1,C_1)
+S_1(\mathbf{a})F_0(V_{0,SOS},C_{0,SOS})\\
-S_0(\mathbf{a})(8-\mathbf{a}^T\mathbf{a}) \mbox{~is SOS},\\
S_0(\mathbf{a})\mbox{~~is~ SOS}.
\end{array}
\right.
\end{array}
\eeas
To overcome the non-convexity due to the term $u_1^2$ in $F_1$ is issue, we use the Schur complement formula, i.e., $F_1\le 0$ is equivalent to the non-negativeness of the matrix
\beas
E_1\eqdef
\left[
\begin{array}{cc}
-W_1&u_1\\
u_1&1
\end{array}
\right]
\eeas
where $W_1=\mathbf{f}\cdot\nabla_{\mathbf{a}}V_1+( \mathbf{e} u_1+F\mathbf{a} u_1+\mathbf{m}\frac{du_1}{d\mathbf{a}}\mathbf{f} )\cdot\nabla_{\mathbf{a}}V_0 - C_1$.
Thus, the problem $O_1'$ can be transformed to
\beas
&O_1'': \begin{array}{c}
\displaystyle \min_{V_1,u_1,S_0,S_1, C_1} C_1, \mbox{~~s.t.~~} \\
[2ex]
\left\lbrace
\begin{array}{c}
\bfz^TE_1(V_{0,SOS},V_1,u_1,C_1)\mathbf{z}
+S_1(\mathbf{a},\mathbf{z})F_0(V_{0,SOS},C_{0,SOS})\\
-S_0(\mathbf{a},\mathbf{z})(8-\mathbf{a}^T\mathbf{a}) \mbox{~is SOS},
\\
S_0(\mathbf{a},\mathbf{z})\mbox{~~is~ SOS}.
\end{array}
\right.
\end{array}
\eeas
To reduce the complexity of the optimization, we choose $S_0(\mathbf{a},\mathbf{z})=s_0(\mathbf{a})\mathbf{z}^T\mathbf{z}$ and $S_1(\mathbf{a},\mathbf{z})=s_1(\mathbf{a})\mathbf{z}^T\mathbf{z}$.

When $d_{u_1}=1$, $d_{s_1}=2$ and $d_{s_0}=d_{V_1}=4$, solving $O_1''$ yields $C_{1,SOS}=-15.5$, associated with a linear controller
\begin{equation}
\begin{array}{ll}
u_{1,SOS}=&-0.0845 x _{1}+1.6704 x _{2}+0.9074 x _{3}+0.5520 x _{4}\\
&-2.8371 x _{5}+1.5102 x _{6}-2.2499 x _{7}+2.9690 x _{8}\\
&-0.1860 x _{9}-3.5105 x _{10}
\end{array}
\label{linearcontrol}
\end{equation}
Further, solving $O_1''$ for $d_{u_1}=2$ results in a quadratic controller with $C_{1,SOS}=-18.49$ and
\begin{equation}
u_{1,SOS}=\bfk\mathbf{a}+\mathbf{a}^TM\mathbf{a}.
\label{quadraticcontrol}
\end{equation}
The numerical values of $\bfk$ and $M$ are given in Appendix C. 
The quadratic controller results in a lower value of $C_{1, SOS},$ as it could be expected.

\subsection{Closed-loop control results for the reduced-order model}
Once the state-feedback controller coefficients were identified, the actual time-averaged costs $\overline{\Phi}$ and $\overline{\Phi}_0$ were calculated from numerical integration of the ordinary differential equation of the closed-loop ROM (\ref{eq:ode-sys}), for a total integration time long enough to provide converged values of the averages and for several increasing values of the parameter $\epsilon$. A point on the periodic orbit of the ROM was selected as initial condition and initial transients were discarded from the calculation of the time average to improve the accuracy.

\begin{figure}[!h]
  \centerline{\includegraphics[width=0.75\textwidth]{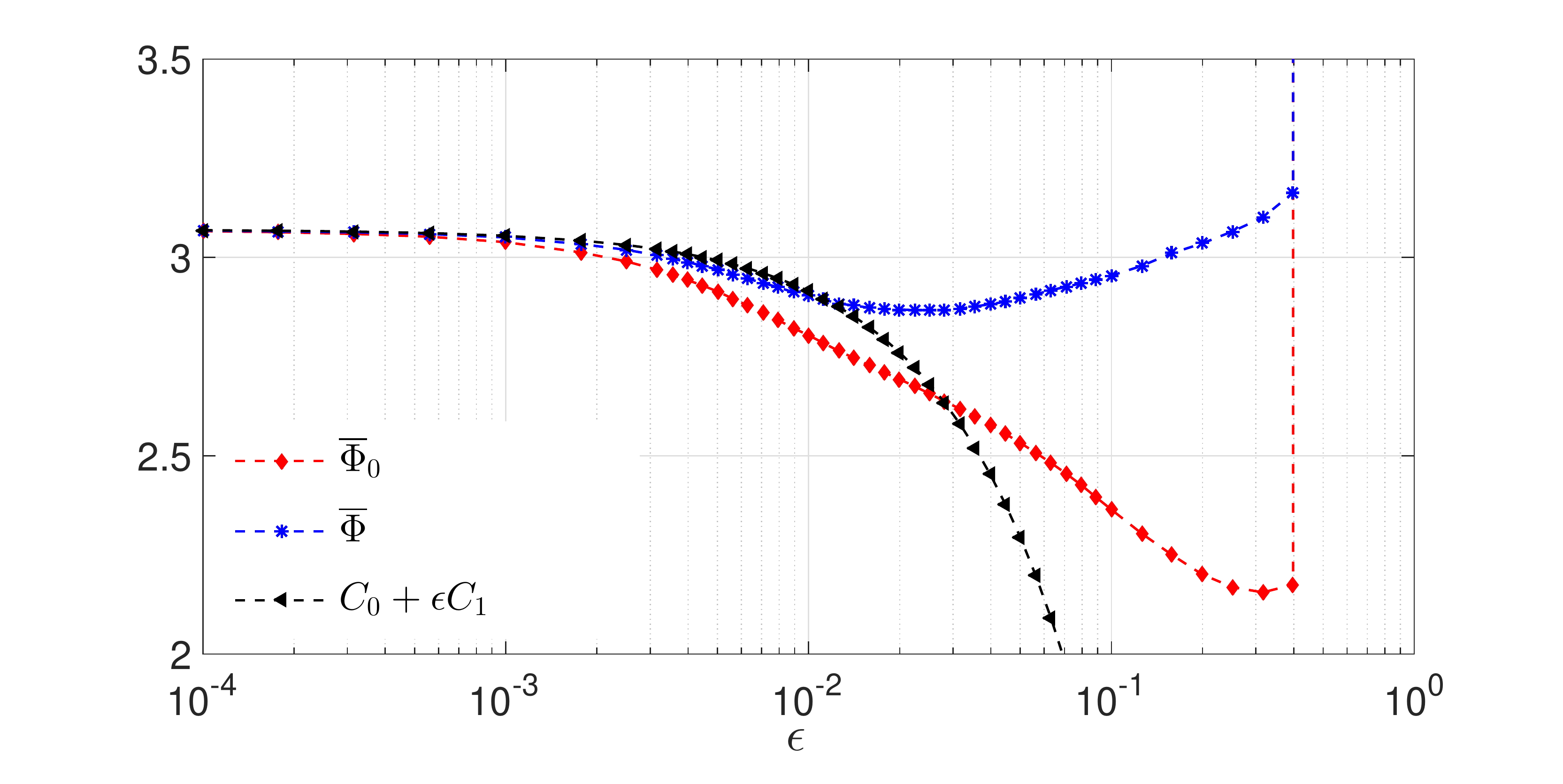}}
  \caption{Time-averaged costs $\overline{\Phi}_0$ and $\overline{\Phi}$ from closed-loop ROM simulation  using linear controllers, and the asymptotics $C_{0,SOS}+\epsilon C_{1,SOS}$ of the upper bound for  $\overline{\Phi}.$}
\label{fig:11}
\end{figure}

\begin{figure}
  \centerline{\includegraphics[width=0.75\textwidth]{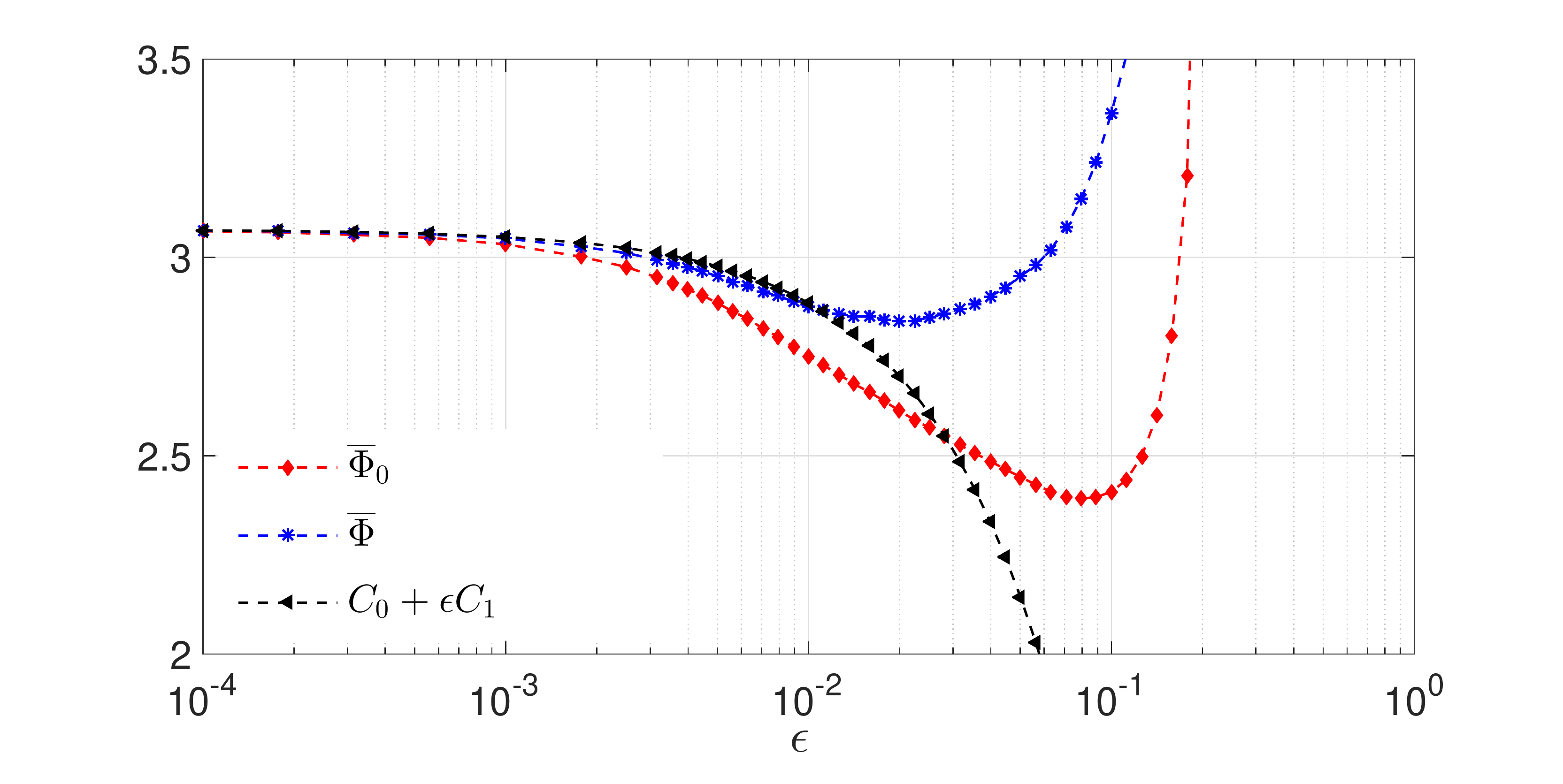}}
\caption{Time-averaged costs $\overline{\Phi}_0$ and $\overline{\Phi}$ from closed-loop ROM simulation  using quadratic controllers, and the asymptotics $C_{0,SOS}+\epsilon C_{1,SOS}$ of the upper bound for  $\overline{\Phi}.$}
\label{fig:12}
\end{figure}

Results are summarized in Figs.~\ref{fig:11} and \ref{fig:12}, for the linear and quadratic controllers (\ref{linearcontrol}) and (\ref{quadraticcontrol}) respectively. The linear approximation of the bound $C_{SOS}$, i.e., $C_{0,SOS}+\epsilon C_{1,SOS},$ is also presented for comparison.

Numerical simulation of the closed-loop system shows that the linear small-feedback controller $u=\epsilon u_{1,SOS}$ reduces the actual long-time average cost to a minimum value equal to $\overline{\Phi} = 2.867$, for $\epsilon = 0.025$. At this value, the system's energy $\overline{\Phi}_0$ is equal to 2.656. The system's energy decreases remarkably more if the parameter $\epsilon$ is increased further, as the control effort becomes significant, as hinted by the increase of the total cost $\overline{\Phi}$. We point out that the difference between $\overline{\Phi}$ and $\overline{\Phi}_0$ is the time average of $u^2/\epsilon$. Hence the penalisation factor $1/\epsilon$ decreases monotonically as the magnitude of the control is increased, whereas one might be interested in having a constant penalisation factor that has a physical, rather than technical, meaning.

The two-term expansion $C_{0,SOS}+\epsilon C_{1,SOS}$ is only a linear approximation of $C_{SOS}$. Thus, formally, it behaves correctly as an upper bound of $\overline{\Phi}$ only when $\epsilon$ is small, but the approximation breaks down when $\epsilon$ is further increased.
 In fact, the results of the numerical calculations suggest that  at $\epsilon=0$ the slope $d\overline{\Phi}(\epsilon)/d\epsilon <C_{1, SOS} = -15.5$.

The effect of $\epsilon$ on $\overline{\Phi}$ can be seen more clearly by investigating the qualitative properties of the long-term behaviour of the closed-loop system. For $\epsilon = 0$ the trajectories of the system converge to a stable limit cycle, over which the long-time average of $\Phi = \Phi_0$ is 3.07. When $\epsilon$ is increased the controller reduces the ``size'' of this limit cycle, where the ``size'' is measured by $\Phi_0$. The cost is reduced up to $\epsilon\approx0.4$, after which the cost grows dramatically as the controller induces an internal bifurcation in the system.

The quadratic small-feedback controller, Fig.~\ref{fig:12}, yields qualitatively similar results. Although the slope of $\overline{\Phi}(\epsilon)$ at $\epsilon=0$ is larger than that associated with  linear controllers, the minimum value of $\overline{\Phi}$ for a finite $\epsilon = 0.02$ is 2.838 at which $\overline{\Phi}_0 = 0.2613.$ This minimal value of~$\overline{\Phi}$ is only marginally smaller than the minimal value obtained with a linear controller. For $\epsilon > 0.08$ the cost increases significantly.

\begin{figure}
  \centerline{\includegraphics[width=0.5\textwidth]{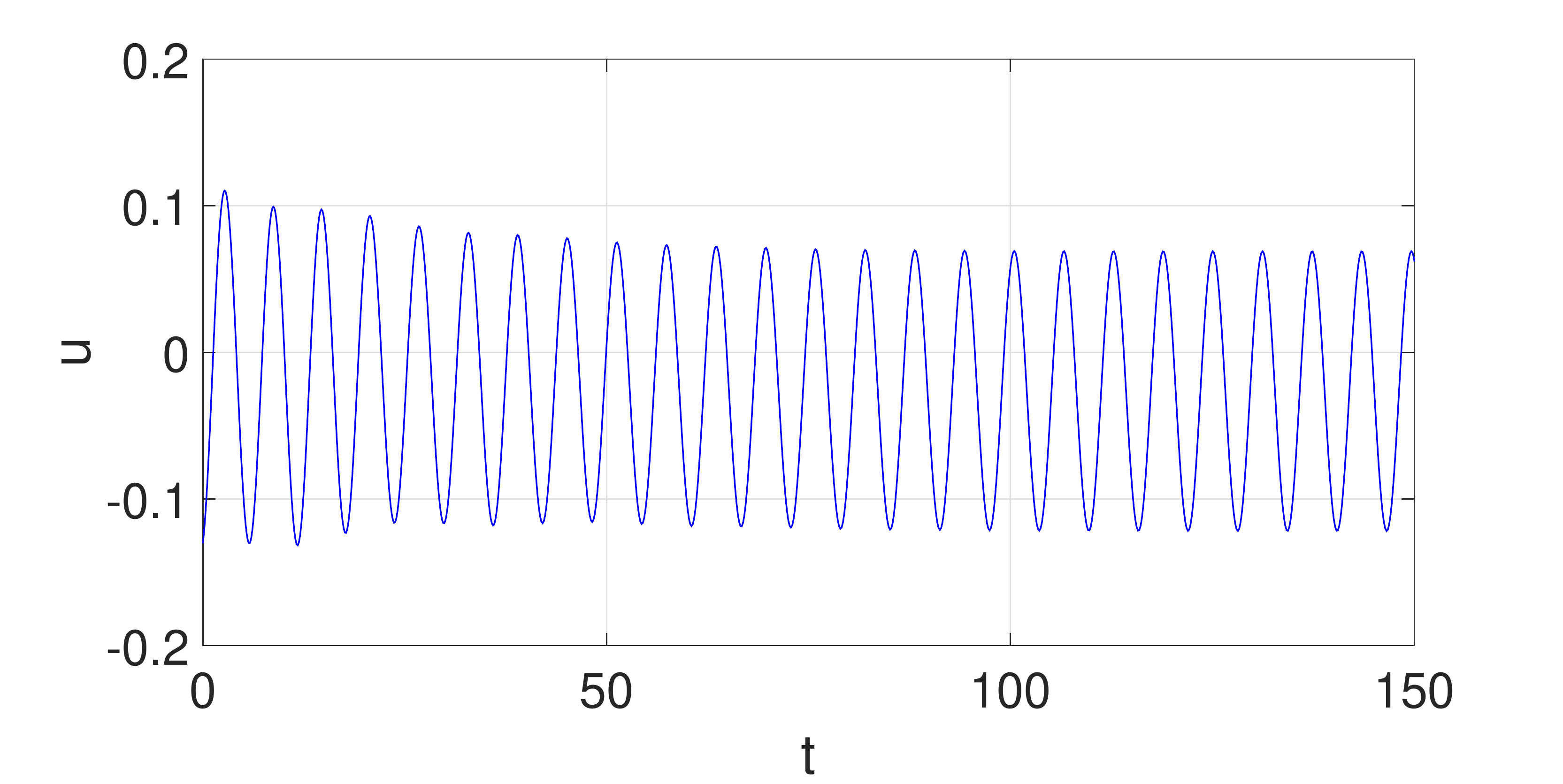}}\hfill
\caption{Profile of the control input for the linear feedback controller with $\epsilon=0.025$.}
\label{fig:21}
\end{figure}

\begin{figure}
   \centerline{\includegraphics[width=0.5\textwidth]{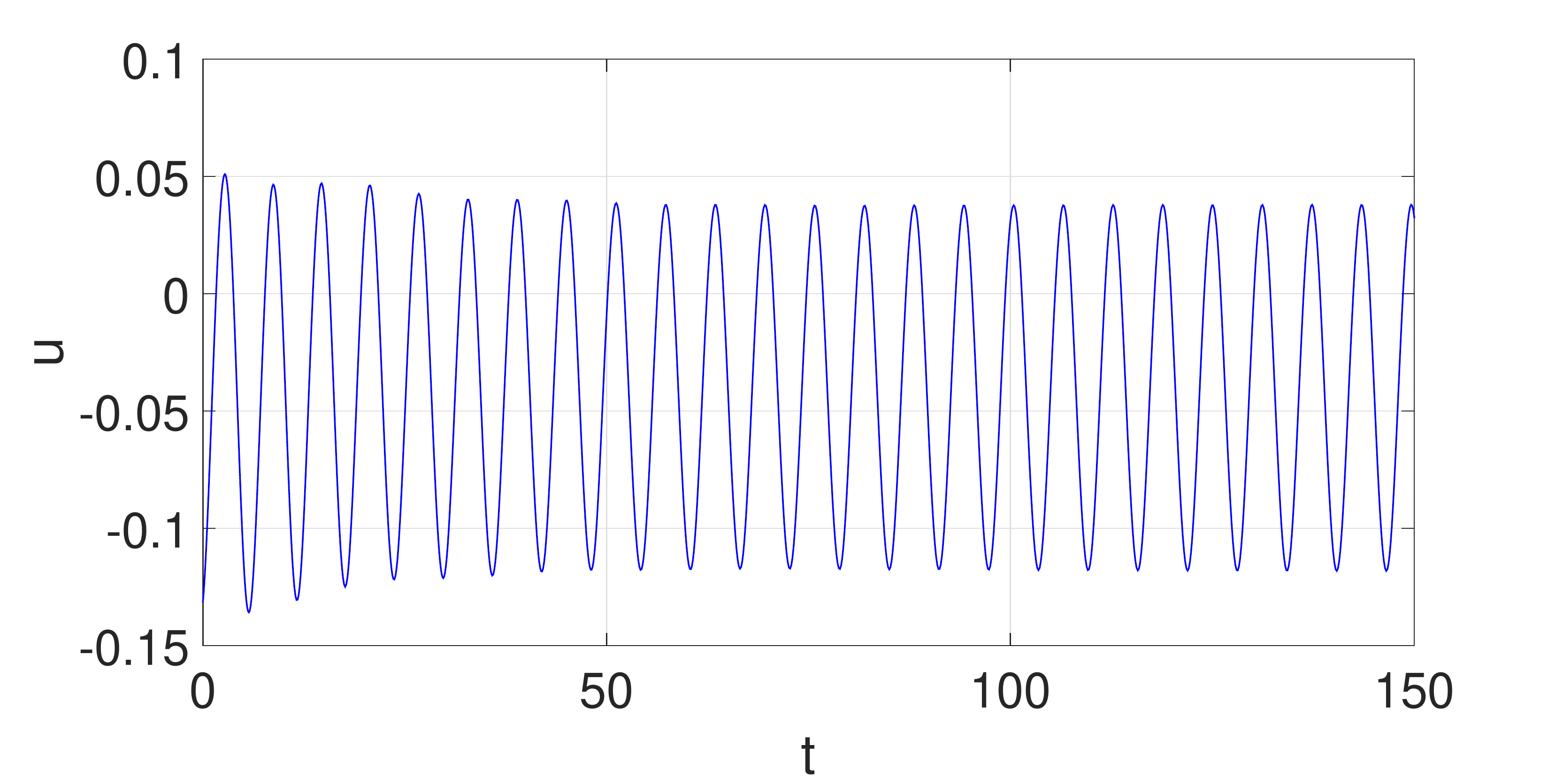}}
\caption{Profile of the control input for the quadratic feedback controller with $\epsilon=0.02$.}
\label{fig:22}
\end{figure}

Figures~\ref{fig:21} and \ref{fig:22} show the control input profiles when $\epsilon=0.025$ in the linear case and $\epsilon=0.02$ in the quadratic case, respectively, the values at which the largest reduction of the long-time averaged cost was obtained in numerical integration of the closed-loop ROM.
In both situations, the control input is a periodic function of time, with a period of oscillation approximately equal to six nondimensional time units. This occurs because when the system trajectory is on the limit cycle, the state variables are periodic, and so is the control input. This also applies when the system is controlled, and its trajectory is on a different, controlled limit cycle. Interestingly, the effect of the quadratic terms in the quadratic controller, which would introduce frequency components in $u$ at twice the fundamental frequency, appears small, possibly justifying the similarity in performance between linear and quadratic controllers.

In summary, the proposed small-feedback controller designed to yield a reduced bound of the long-time averaged cost for small $\epsilon$ also reduces the long-time average cost itself. 

\commentout{ Note that the dThe efficiency of the proposed control scheme was also revealed by the following two facts: (1) the complex dynamics of the flow model makes the control task difficult, and (2) not only linear feedback but also nonlinear feedback can be considered in the same way.}

\subsection{Closed-loop control results for direct numerical simulation}
The linear and quadratic controllers (\ref{linearcontrol}) and (\ref{quadraticcontrol}) were also implemented in direct numerical simulation of the fluid flow, to assess their performance. Multiple simulations, for increasing values of $\epsilon$, have been run similarly to the results shown in Figs.~\ref{fig:11} and~\ref{fig:12}.  The initial condition is chosen to lie on the periodic limit cycle of the full-order system associated with vortex shedding. A total integration time sufficient to have converged long-time averages, discarding the initial transient, was selected.

\begin{figure}
	\centering
	\includegraphics[width=0.5\textwidth]{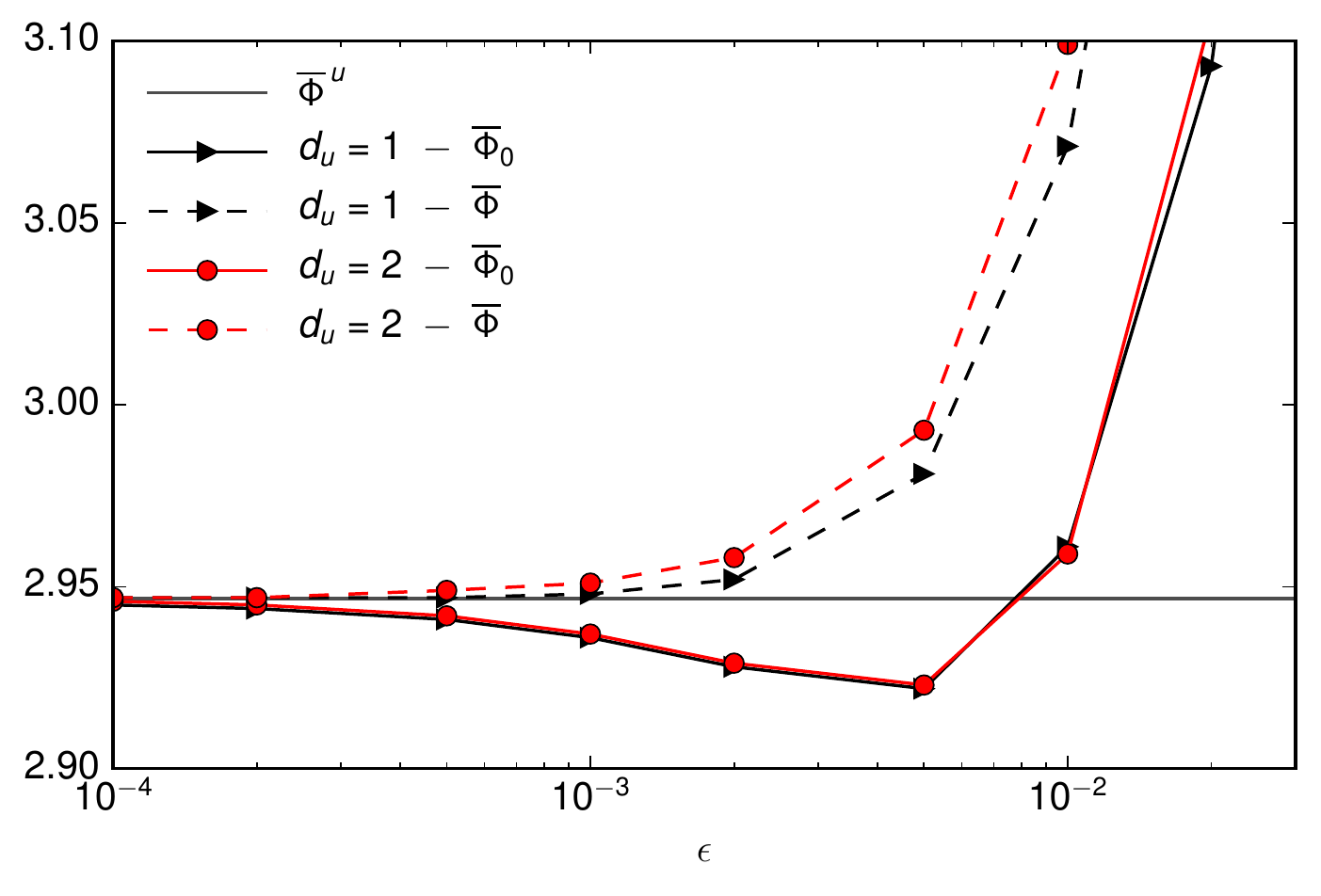}
	\caption{Long-time averaged cost as a function of $\epsilon$ obtained by closed-loop direct numerical simulation. Results for linear (triangles) and quadratic (circles) controllers are reported. }
	\label{fig:dns-performance}
\end{figure}

Figure \ref{fig:dns-performance} summarizes the results for linear and quadratic controllers in direct numerical simulation, and it is the direct equivalent of Figs.~\ref{fig:11} and \ref{fig:12}.
The solid curves refer to $\overline{\Phi}_0$, whereas the dashed curves refer to $\overline{\Phi}$, hence including the control penalization $\overline{u^2}/\epsilon$. The horizontal line for $\overline{\Phi}^{\,u}$ is the value of the long-time averaged cost for the uncontrolled system.

It can be observed that the time-averaged system's energy $\overline{\Phi}_0$ initially decreases as $\epsilon$ increases, up to $\epsilon=5\times10^{-3}$. However, if the parameter $\epsilon$ is further increased, performance worsens. This value is lower by a factor of $4\sim 5$ than what found in simulation of the ROM, in Figs.~\ref{fig:11} and \ref{fig:12}. As a result, the minimum value of $\overline{\Phi}_0$ is 2.922, and the percentage reduction of the long-time averaged cost is lower than that observed from application of the same controller on the ROM.
Note that without control, the time-averaged resolved energy associated with the attractor of the full-order system, 2.947, is slightly lower than that associated with the ROM, 3.07, because the effects of truncation of high-order POD modes have not been completely recovered by the calibration. Also note that the slope of $\Phi_0$ as a function of $\epsilon$ evaluated at $\epsilon=0$, (estimated from a linear fit of the first three data points), is -10.0, for both the linear and quadratic controllers. The value obtained from closed-loop simulation of the ROM, deduced using a similar method, is -33.134 for the linear controllers, and -38.322 for the quadratic controllers, hence significantly lower.

If the cost of the control is also included, no reduction of the total cost is observed in DNS. However, it is worth pointing out, as anticipated, that the formalism introduced to force the smallness of the control, with the penalisation on the input as $u^2/\epsilon$, is artificial.  A further important result is that, in DNS, the effect of quadratic controllers is similar to that of linear controllers, but require a higher control input and result in a higher total cost. However, for the ROM, quadratic controllers resulted in slightly better performance. This result seems to suggest that high-degree polynomial controllers designed on approximate models might not have better performances in DNS.

\begin{figure}[h]
	\centering
	\includegraphics[width=0.75\textwidth]{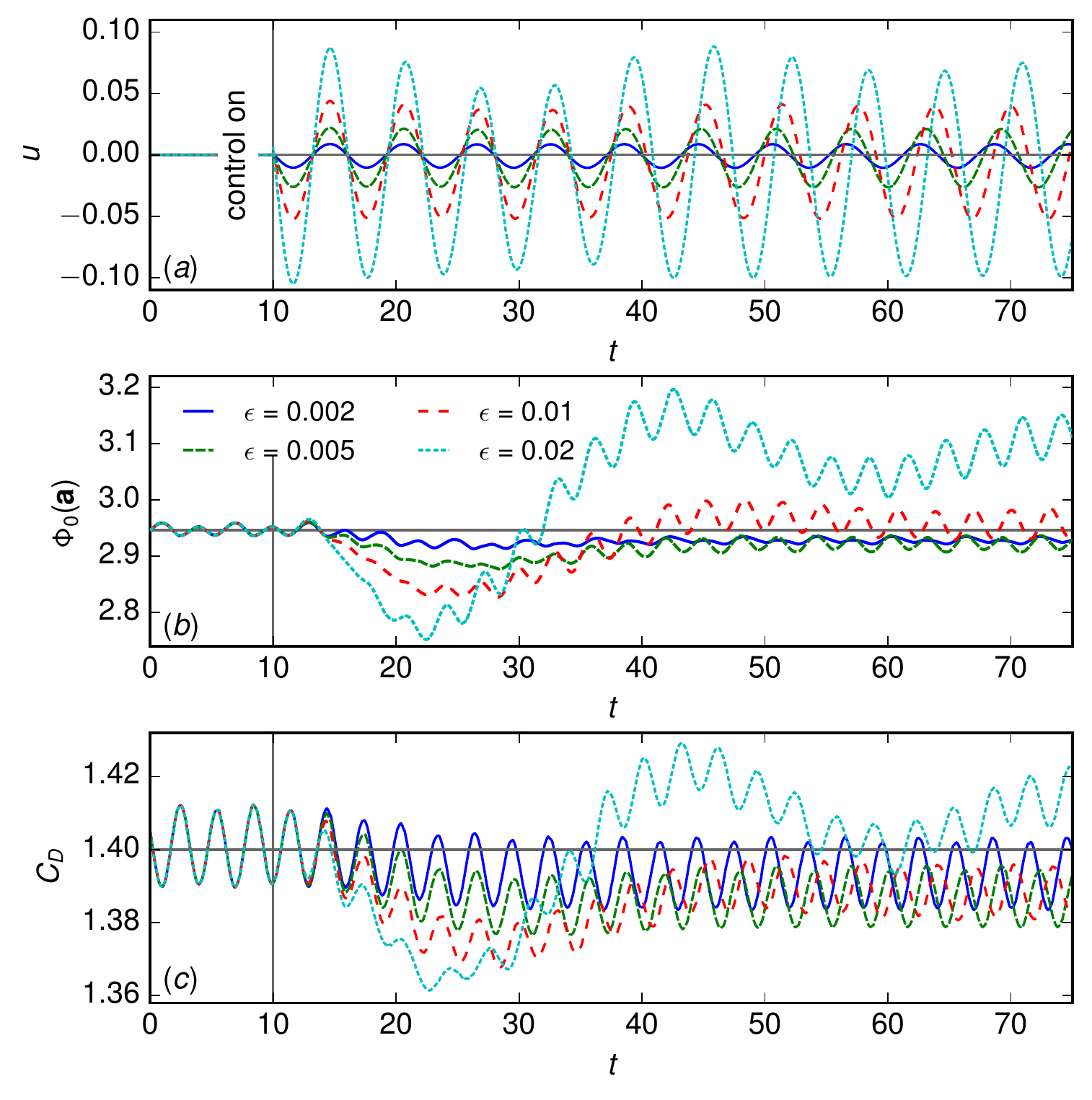}
	\caption{Closed-loop control results for linear controllers with increasing $\epsilon$.
(a): time history of control input, (b): time history of the system's energy, and (c): time history of drag coefficient.}
	\label{fig:control-input-and-cost}
\end{figure}

Figure \ref{fig:control-input-and-cost} shows time histories of the control input $u(t)$, panel $(a)$, of the the system's energy $\Phi_0(\mathbf{a})$, panel $(b)$, and of the drag coefficient $C_D(t)$, panel $(c)$, obtained from  direct numerical simulation results of linear controllers for four different values of $\epsilon$ selected from Fig.~\ref{fig:dns-performance}. Control is activated at $t=10$. The drag coefficient is obtained from integration of the pressure and shear stress distributions along the cylinder surface.

The input $u(t)$ is a periodic function of time. For \mbox{$\epsilon=5\times 10^{3}$}, for which the long-time average of $\Phi_0$ is reduced the most, its peak-to-peak variation is small, on the order of 0.05, i.e., the cylinder oscillates only slightly. Time integration of $u(t)$, yielding the cylinder angular rotation, shows that the total oscillation is about 26 degrees for this case. The control input is periodic because it is a linear function of the state variables, some of which are close to sinusoids where the system trajectory lies on the periodic limit cycle associated with vortex shedding. This applies when the control is small and the dynamics are not perturbed significantly. For larger values of $\epsilon$, a modulation of the control signal is observed as the actuated dynamics of the full-order system change significantly.

When control is activated the system exhibits a transient during which $\Phi_0$ and $C_D$ initially decrease. Their reduction and the time rate at which they decrease is considerable for larger $\epsilon$. This result clearly shows that the small-feedback controller has correctly identified the right physical mechanism for control to mitigate vortex shedding when the system is near the attractor.

However, after approximately $10\sim 20$ time units, approximately $2\sim 3$ shedding cycles, when the system trajectory departs from the neighbourhood of the attractor, control effectiveness is lost and the cost begins to increase, especially for larger values of $\epsilon$. For small values, i.e. small control input, the cost eventually settles down to a time-averaged value lower than that associated with the uncontrolled system, indicated in panels $(b)$ and $(c)$ as a horizontal line. For larger values of $\epsilon$, e.g. 0.02, the cost increases significantly. This occurs because the small-feedback controller is designed to control the flow only in a narrow volume of the phase space containing the attractor. For large $\epsilon$ the controller drives the fluid flow too far from the model design point and performance worsen. From this perspective, the relatively good reduction observed in Figs.~\ref{fig:11}  and \ref{fig:12} for large $\epsilon$ is artificial as the ROM dynamics for such a large control input are not physically realistic.

The fundamental origin of the discrepancy in performance between the ROM and the full-order system is that the present ten-mode reduced-order model is constructed to model the dynamics of the actuated flow only in a relatively small volume of the phase space around the design point, i.e. around the unactuated periodic limit cycle. Although the effect of actuation are in part taken into account in the construction of the model, and in particular in the generation of the POD subspace, when the amplitude of the control is increased via $\epsilon$ the separation between actual and modelled dynamics increases consistently. A larger reduced-order model would give accurate prediction in a larger volume of the phase space, that is for larger values of $\epsilon$.  If this were the case, a better match of the actuated dynamics near the attractor, i.e. a better match between the results in Figs.~\ref{fig:11}, \ref{fig:12} and \ref{fig:dns-performance} would have been observed.  Hence, the relatively poor performance observed in DNS is almost entirely due to the necessity to include only ten modes in the Galerkin projection. We have derived reduced-order models for this flow problem using 20 or 30 modes that could describe with more accuracy the actuated dynamics of the flow. However, the computational costs associated with the solution of the SOS problems $O_0',\ O_1''$, memory- and time-wise, using these larger models was too large. Further algorithmic and numerical advances are required to enable the application of the present methodology to larger systems.

\section{Conclusion}\label{seq:Conclusion}

Based on sum-of-squares decomposition of polynomials and semidefinite programming, a numerically-tractable approach is presented for expensive control of the  long-time average cost of polynomial dynamical systems. The control law is restricted to polynomials of the system state. The derivation of the controller is given in terms of solvability conditions of state-dependent linear and bilinear inequalities. The non-convexity in SOS optimization is resolved by making use of the smallness of the perturbation parameter describing the (high) control cost without applying iterative algorithms.
The proposed control design scheme has been applied to the problem of mitigating two-dimensional vortex shedding past a circular cylinder at low Reynolds number, using controlled rotary oscillations in a full-information controller setting. Linear and nonlinear controllers, designed using an approximate reduced-order model of the actuated flow, have been also implemented in direct numerical simulation.

The proof of concept of the idea of using the upper bound of long-time average cost control as the objective of the control design, the method of overcoming the non-convexity of simultaneous optimization of the control law and the tunable function under the structure of expensive control, as well as the detailed implementation of the proposed control scheme to a particular fluid flow are the three main contributions of the present paper.

\section*{Acknowledgment}
Funding from EPSRC under the grants EP/J011126/1 and EP/J010073/1 and support
in kind from Airbus Operation Ltd., ETH Zurich (Automatic Control Laboratory), University of Michigan (Department of Mathematics), and
University of California, Santa Barbara (Department of Mechanical Engineering) are gratefully acknowledged.

\section*{Appendix A: Proof of Theorem 1}

\noindent Let $V_{SOS}=V_{0,SOS}+\epsilon V_{1,SOS}$.
By substituting $V=V_{SOS}, C=C_{\kappa, SOS}, \bfu=\bfu_{SOS}$ in the constraint function $F(V,\bfu,C)$ that is defined in (\ref{cc1}), the remaining task
is to seek small $\epsilon>0$ such that
\bea
F(V_{SOS},\bfu_{SOS},C_{\kappa, SOS})
\le
0.
\label{cc1_1}
\eea
Notice that
\bea
\begin{array}{rl}
&F(V_{SOS},\bfu_{SOS},C_{\kappa,SOS})\\
&\ \ =
F_0(V_{0,SOS},C_{0,SOS})
+\epsilon F_1(V_{0,SOS},V_{1,SOS},\bfu_{1,SOS},C_{1,SOS})\\
&\ \ \ \ \ +\epsilon (1-\kappa)C_{1,SOS}+\epsilon^2 w(\bfx,\epsilon)
\end{array}
\label{boundx}
\eea
where
\beas
\begin{array}{ll}
w(\bfx,\epsilon)=&\bfG\bfu_1\cdot  \nabla_{\bfx}V_{1,SOS}\\
&+\frac{1}{\epsilon^2}\left(\Phi(\bfx,\epsilon \bfu_{1,SOS})-\Phi(\bfx,0)-\epsilon\frac{\partial \Phi}{\partial \bfu}(\bfx,0)\bfu_{1,SOS}\right)
\end{array}
\eeas
and $F_0, F_1$, being polynomial in $\bfx$, possess all the continuity properties implied by the proof.
Let $\mathcal{D}\subseteq{\mathbb R}^n$ be the phase domain that interests us, where the closed-loop trajectories are all bounded.
Then,
\bea
F_{1,max}\eqdef \max_{\bfx\in\mathcal{D}}F_1(V_{0,SOS},V_{1,SOS},\bfu_{1,SOS},C_{1,SOS})<\infty,
\label{bound1}
\eea
and $w(\bfx,\epsilon)$ is bounded for any $\bfx\in\mathcal{D}$ and any finite $\epsilon$ (the latter following from the standard mean-value-theorem-based formula for the Lagrange remainder). By (\ref{boundx}) and (\ref{bound1}),
\bea
\begin{array}{ll}
F(V_{SOS},\bfu_{SOS},C_{\kappa,SOS})\le
&F_0(V_{0,SOS},C_{0,SOS})+\epsilon F_{1,max}\\&+\epsilon (1-\kappa)C_{1,SOS}+O(\epsilon^2)
\end{array}
\label{bound1_ex}
\eea
Consider the two inequality constraints obtained by solving $O_0$ and $O_1$:
\bea
\left\{
\begin{array}{c}
 F_0(V_{0,SOS},C_{0,SOS})\le 0, \\
[1ex]
 F_1(V_{0,SOS},V_{1,SOS}, \bfu_{1,SOS},C_{1,SOS})\le 0\\
~~\forall \bfx ~~\mbox{such that}~ F_0(V_{0,SOS},C_{0,SOS})= 0.
\end{array}
\right.
\label{cons}
\eea

Define $\mathcal{D}_{\delta}\eqdef \left\{\bfx\in \mathcal{D} ~|~ \delta\le F_0(V_{0,SOS},C_{0,SOS})\le 0\right\}$ for a given constant $\delta\le 0$.
Clearly, $\mathcal{D}_{\delta}\rightarrow \mathcal{D}_{0}$ as $\delta\rightarrow 0$. Further define
\[
F_{1,\delta}(\delta)\eqdef \max_{\bfx\in \mathcal{D}_{\delta}} F_1(V_{0,SOS},V_{1,SOS},\bfu_{1,SOS},C_{1,SOS}).
\]

By the second constraint in (\ref{cons}), $\lim_{\delta\rightarrow 0}F_{1,\delta}(\delta)\le 0$. Therefore, by continuity and the fact $C_{1,SOS}<0$, for any $0<\kappa<1$ there exists a constant $\delta_{\kappa}<0$ such that
\bea
\begin{array}{r}
F_1(V_{0,SOS},V_{1,SOS}, \bfu_{1,SOS},C_{1,SOS})\le F_{1,\delta_{\kappa}}<-\frac{1}{2}(1-\kappa)C_{1,SOS}\\
\forall \bfx\in \mathcal{D}_{\delta_{\kappa}}.
\end{array}
\label{bound3}
\eea
In consequence, (\ref{boundx}), the first constraint in (\ref{cons}), and (\ref{bound3}) render to
\bea
\begin{array}{rl}
& F(V_{SOS},\bfu_{SOS},C_{\kappa,SOS}) \\
 \le&F_0(V_{0,SOS},C_{0,SOS})+\epsilon F_{1,\delta_{\kappa}}+\epsilon (1-\kappa)C_{1,SOS}+O(\epsilon^2) \label{bound4'}\\
\le&\frac{\epsilon}{2} (1-\kappa)C_{1,SOS}+O(\epsilon^2)\le 0, ~~\forall \bfx\in \mathcal{D}_{\delta_{\kappa}}
\end{array}
\eea
for sufficiently small $\epsilon$. Hence, (\ref{cc1_1}) holds for any $\bfx\in \mathcal{D}_{\delta_{\kappa}}$.

The remaining is to prove (\ref{cc1_1}) for any $\bfx\in \mathcal{D}\setminus\mathcal{D}_{\delta_{\kappa}}$. By the definition of the set $\mathcal{D}_{\delta_{\kappa}}$, we have
\bea
F_0(V_{0,SOS},C_{0,SOS})<\delta_{\kappa}<0,  ~~\forall \bfx\in \mathcal{D}\setminus\mathcal{D}_{\delta_{\kappa}}.
\label{bound5}
\eea
Then, (\ref{bound1_ex}) and (\ref{bound5}) yield
\bea
\begin{array}{rl}
F(V_{SOS},\bfu_{SOS},C_{\kappa,SOS})\le&
\delta_{\kappa}+\epsilon F_{1,max}+\epsilon (1-\kappa)C_{1,SOS}\\&+O(\epsilon^2)\\
\le &\delta_{\kappa}+O(\epsilon)\le 0, ~~\forall \bfx\in \mathcal{D}\setminus\mathcal{D}_{\delta_{\kappa}}
\end{array}
\label{bound1_ex1}
\eea
if $\epsilon$ is sufficiently small.Then,  (\ref{bound4'}) and (\ref{bound1_ex1}) imply that (\ref{cc1_1}) holds $\forall \bfx\in \mathcal{D}$. The proof is completed.

\section*{Appendix B: Galerkin projection terms}
\label{sec:app-ode-sys}

\noindent
With $\omega_i(\mathbf{a})$ being the scalar vorticity field associated with mode $\mathbf{v}_i(\mathbf{a})$, and similarly for $\overline{\mathbf{v}}$ and $\mathbf{v}_c(\mathbf{a})$,  Galerkin projection results in the following coefficients:
\[
	c_i = -\frac{1}{Re}\int_\Omega \omega_i \nabla^2 \overline{\omega} \mathrm{d}\Omega - \int_\Omega \mathbf{v}_i \cdot (\overline{\mathbf{v}} \cdot \nabla \overline{\mathbf{v}})\mathrm{d}\Omega,
\]
\[
\begin{array}{ll}
	L_{ij} = &-\frac{1}{Re}\int_\Omega \omega_i \nabla^2 \omega_j \mathrm{d}\Omega - \int_\Omega \mathbf{v}_i \cdot (\overline{\mathbf{v}} \cdot \nabla \mathbf{v}_j)\mathrm{d}\Omega \\
	&- \int_\Omega \mathbf{v}_i \cdot (\mathbf{v}_j \cdot \nabla \overline{\mathbf{v}})\mathrm{d}\Omega,
\end{array}
\]
\[
    N_{ijk} = -\int_\Omega \mathbf{v}_i \cdot (\mathbf{v}_j \cdot \nabla) \mathbf{v}_k \mathrm{d}\Omega,
\]
\[
	m_i = - \int_\Omega \mathbf{v}_i \cdot \mathbf{v}_c \mathrm{d}\Omega,
\]
\[
	e_i = - \int_\Omega \mathbf{v}_i \cdot (\mathbf{v}_c \cdot \nabla \overline{\mathbf{v}} + \overline{\mathbf{v}} \cdot \nabla \mathbf{v}_c ), \mathrm{d}\Omega - \frac{1}{Re} \int_\Omega \omega_i\omega_c  \mathrm{d}\Omega,
\]
\[
	b_i = -\int_\Omega \mathbf{v}_i \cdot (\mathbf{v}_c \cdot \nabla \mathbf{v}_c) \mathrm{d}\Omega,
\]
\[
	F_{ij} = -\int_\Omega \mathbf{v}_i \cdot (\mathbf{v}_j \cdot \nabla \mathbf{v}_c + \mathbf{v}_c \cdot \nabla \mathbf{v}_j) \mathrm{d}\Omega.
\]
In the present case all the coefficients $b_i$ are identically zero because of the radial symmetry of the control function $\mathbf{v}_c$. Domain integrals are evaluated numerically on the triangular unstructured mesh by using a linear approximation of the integrand function based on nodal values. All derivatives are computed using a local quadratic interpolation scheme available in Algorithm 624 from~\cite{renka:04}. Strictly, some of the above definitions do not contain the line integrals on the boundary of the domain arising from the use of vector calculus identities to eliminate the Laplacian, as in Appendix 2 of~\cite{be:05}, as these are found to be quite small and negligible in the present case with respect to the domain integrals above. Appropriate symmetries in the tensor $N_{ijk}$ are numerically enforced after the computations of the integrals to ensure that the nonlinear term is energy preserving, (see e.g.~\cite{schlegel:15} for a discussion on this topic for the present case).

\section*{Appendix C: quadratic controller}

The numerical values of of $\bfk$ and $M$ in (\ref{quadraticcontrol}) are:

\begin{equation}\label{K_matrix}
\bfk=\left[
\begin{array}{ >{\scriptstyle   \hspace{-4pt}} r  >{\scriptstyle   \hspace{-8pt}} r  >{\scriptstyle   \hspace{-8pt}} r  >{\scriptstyle   \hspace{-8pt}} r  >{\scriptstyle   \hspace{-8pt}} r  >{\scriptstyle   \hspace{-8pt}} r  >{\scriptstyle   \hspace{-8pt}} r  >{\scriptstyle   \hspace{-8pt}} r  >{\scriptstyle   \hspace{-8pt}} r  >{\scriptstyle   \hspace{-8pt}} r }
   -0.0790&
    1.5406&
    1.0863&
    0.5558&
   -2.4814&
    0.9818&
   -2.4206&
   -0.1333&
    0.3233&
    0.3959
\end{array}\right]
\end{equation}

\begin{equation}\label{M_matrix}
M=\scriptstyle
\left[\scriptstyle
\begin{array}{ >{\scriptstyle   \hspace{-4pt}} r  >{\scriptstyle   \hspace{-8pt}} r  >{\scriptstyle   \hspace{-8pt}} r  >{\scriptstyle   \hspace{-8pt}} r  >{\scriptstyle   \hspace{-8pt}} r  >{\scriptstyle   \hspace{-8pt}} r  >{\scriptstyle   \hspace{-8pt}} r  >{\scriptstyle   \hspace{-8pt}} r  >{\scriptstyle   \hspace{-8pt}} r  >{\scriptstyle   \hspace{-8pt}} r }
  -0.2874 &   0.0285 &  -0.0194 &  -0.1038  & -0.0775 &   0.0181 &  -0.1399  &  0.0369 &   0.0854  & -0.0277\\
    0.0285 &  -0.2981&   -0.0784 &   0.0302 &  -0.0186&   -0.1099 &   0.0560 &   0.0630 &  -0.1956 &   0.2916\\
   -0.0194  & -0.0784&   -0.2135 &   0.0097  & -0.0573 &   0.0516 &  -0.0578 &   0.1121 &  -0.6484 &   0.4140\\
   -0.1038  &  0.0302 &   0.0097 &  -0.2520 &   0.0537 &   0.0476 &  -0.1767 &  -0.0604 &  -0.0824 &   0.0619\\
   -0.0775  & -0.0186 &  -0.0573 &   0.0537 &  -0.2072 &  -0.0730 &  -0.3547  & -0.0262 &   0.4199 &   0.2580\\
    0.0181 &  -0.1099 &   0.0516 &   0.0476  & -0.0730 &  -0.1572 &   0.2190 &   0.3750 &   0.3342 &  -0.1550\\
   -0.1399  &  0.0560 &  -0.0578 &  -0.1767 &  -0.3547 &   0.2190  & -0.5368  & -0.2473 &  -0.6232 &   0.2140\\
    0.0369  &  0.0630 &   0.1121  & -0.0604  & -0.0262 &   0.3750 &  -0.2473 &   0.6261 &  -0.2342 &   0.2815\\
    0.0854  & -0.1956 &  -0.6484&   -0.0824 &  0.4199  &  0.3342  & -0.6232  & -0.2342  &  0.2063  &  0.0571\\
   -0.0277  &  0.2916 &   0.4140  &  0.0619  &  0.2580 &  -0.1550 &   0.2140 &   0.2815 &   0.0571  &  0.0221
\end{array}\right]
\end{equation}

\newcommand\Journal{J.~}
\newcommand\International{Int.~}

\end{document}

%% file: Eedefs.tex
\newcommand{\ba}{\begin{array}}
\newcommand{\ea}{\end{array}}
\newcommand{\babc}{\begin{abc}}
\newcommand{\eabc}{\end{abc}}
\newcommand{\bc}{\begin{center}}
\newcommand{\ec}{\end{center}}
\newcommand{\be}{\begin{equation}}
\newcommand{\ee}{\end{equation}}
\newcommand{\bea}{\begin{eqnarray}}
\newcommand{\eea}{\end{eqnarray}}
\newcommand{\beas}{\begin{eqnarray*}}
\newcommand{\eeas}{\end{eqnarray*}}
\newcommand{\bh}{\begin{hangitem}}
\newcommand{\eh}{\end{hangitem}}
\newcommand{\bhi}{\begin{hangitem}}
\newcommand{\ehi}{\end{hangitem}}
\newcommand{\bi}{\begin{itemize}}
\newcommand{\ei}{\end{itemize}}
\newcommand{\bn}{\begin{enumerate}}
\newcommand{\en}{\end{enumerate}}
\newcommand{\bq}{\begin{quote}}
\newcommand{\eq}{\end{quote}}
\newcommand{\btb}{\begin{tabular}}
\newcommand{\etb}{\end{tabular}}
\newcommand{\eqdef}{\stackrel{\triangle}{=}} 
\def\litem[#1]{\item[#1\hfill]}         
\newtheorem{theorem}{Theorem}

\newtheorem{lemma}{Lemma}

\newtheorem{remark}{Remark}

\def\bfa{{\bf a}}

\def\bff{{\bf f}}

\def\bfk{{\bf k}}

\def\bfu{{\bf u}}

\def\bfw{{\bf w}}
\def\bfx{{\bf x}}

\def\bfz{{\bf z}}